\documentclass[10pt]{article}
 \usepackage{color}
\usepackage{amssymb}
\usepackage{amsmath}

\usepackage{amsthm}
\usepackage{epstopdf}
\usepackage{graphicx}
\usepackage{mathtools}
\usepackage[usenames,dvipsnames, table]{xcolor}
\usepackage[hang]{caption}
\usepackage{hyperref}
\usepackage{lscape} 
\usepackage{oldgerm}  
\usepackage{pgfplotstable}
\usepackage{longtable}
\usepackage{array}   

\newcolumntype{L}{>{$}l<{$}} 

\def\be{\begin{equation}}\def\ee{\end{equation}}
\def\bee{\begin{enumerate}}\def\eee{\end{enumerate}}
\def\bei{\begin{itemize}}\def\eei{\end{itemize}}
\def\SU{{\rm SU}}\def\U{{\mathrm U}}\def\SO{{\mathrm SO}}

\newcommand{\nco}{\newcommand}

\nco{\one}{\ensuremath{\,\,\mathrm{l}\!\!\!1}} 
\nco{\NN}{\mathbb{N}}
\nco{\ZZ}{\mathbb{Z}}
\nco{\QQ}{\mathbb{Q}}
\nco{\RR}{\mathbb{R}}
\nco{\CC}{\mathbb{C}}
\nco{\HH}{\mathbb{H}}
\nco{\OO}{\mathbb{O}}

\nco{\red}{\color{red}}
\nco{\blue}{\color{blue}}
\nco{\cyan}{\color{cyan}}
\nco{\brown}{\color{Magenta}}

\nco{\magenta}{\normalcolor}

\nco{\violet}{\color{violet}}
\nco{\redend}{\normalcolor}
\nco{\magentaend}{\normalcolor}

\def\ie{{\it i.e. }}
\def\ommit#1{{}}
\def\({\left(}
\def\){\right)}

\def\ie{{\it i.e.,\/}\ }
\def\ie{{\rm i.e.,\/}\ }
\def\etc{{\rm etc.\/}\ }

\def\be{\begin{equation}}\def\ee{\end{equation}}
\def\bea{\begin{eqnarray}}\def\eea{\end{eqnarray}}

\nco{\rnc}{\renewcommand}
\rnc{\title}[1]{{\Large\bf\mbox{}\\\medskip#1\bigskip\medskip\\}}
\rnc{\author}[1]{{\large #1\smallskip\\}}
\nco{\address}[1]{{\em #1\medskip\\}}

\newcommand{\incircbin}[1]{%
  \mathbin{%
    \mathchoice%
    {\protect\incircint{\displaystyle}{#1}}%
    {\protect\incircint{\textstyle}{#1}}%
    {\protect\incircint{\scriptstyle}{#1}}%
    {\protect\incircint{\scriptscriptstyle}{#1}}%
  }%
}
\newcommand{\incircint}[2]{%
  \ooalign{$#1\bigcirc$\crcr\hidewidth$#1#2$\hidewidth\crcr}%
}
\newcommand{\circleland}{\incircbin{\land}}

\DeclareMathAlphabet{\mathpzc}{OT1}{pzc}{m}{it}
 \def \w{\mathpzc{w}}

\begin{document}

\begin{titlepage}
\begin{center}
\title{About left-invariant geometry\\ and homogeneous pseudo-Riemannian Einstein structures\\ on the Lie group $\SU(3)$}
\medskip
\author{Robert Coquereaux} 
\address{Aix Marseille Univ, Universit\'e de Toulon, CNRS, CPT, Marseille, France\\
Centre de Physique Th\'eorique}

\bigskip\medskip

\today

\begin{abstract}
\noindent {
This is a collection of notes on the properties of left-invariant metrics on the eight-dimensional compact Lie group $\SU(3)$.
Among other topics we  investigate the existence of  invariant  pseudo-Riemannian Einstein metrics on this manifold. 
We recover the known examples (Killing metric and Jensen metric) in the Riemannian case (signature $(8,0)$), as well as a Gibbons et al example of signature $(6,2)$, and we describe a new example, which is Lorentzian (\ie of signature $(7,1)$.
In  the latter case the associated metric is left-invariant, with isometry group $\SU(3) \times \U(1)$, and has positive Einstein constant. It seems to be the first example of a Lorentzian homogeneous Einstein metric on this compact manifold. 

These notes are arranged into a paper that deals with various other subjects unrelated with the quest for Einstein metrics but that may be of independent interest: among other topics we describe the various groups that may arise as isometry groups of left-invariant metrics on $\SU(3)$, 
provide parametrizations for these metrics, give several explicit results about the curvatures of the corresponding Levi-Civita connections, discuss modified Casimir operators (quadratic, but also cubic) and Laplace-Beltrami operators. 
In particular we discuss the spectrum of the Laplacian for metrics that are invariant under $\SU(3)\times \U(2)$, a subject that may be of interest in particle physics.
}
\end{abstract}
\end{center}

\vspace*{70mm}
\end{titlepage}


\section{Introduction}

This paper is an excursion in the land of left-invariant geometry, more precisely in the land of left-invariant pseudo-Riemannian metrics on the Lie group $\SU(3)$.
Its main purpose is to illustrate several known concepts and methods of Riemannian geometry in this particular case.

The only mathematical result which is probably new is the existence and description of an homogeneous Einstein metric (actually a family), with Lorentz signature, \ie with signature $(7,1)$.
Of course, any quadratic form of signature $(7,1)$ in $R^8$ gives rise, by group translations, to a Lorentzian homogeneous metric on the  group $\SU(3)$, such metrics, usually, are  not Einstein metrics.
The example that we present in sect.~\ref{sec: EinsteinMetrics} seems to be the first example of an homogeneous Einstein Lorentzian metric on this $8$-dimensional compact manifold.
Its isometry group is isomorphic with $\SU(3)\times \U(1)$.

It is probably the proper place to mention that, for us, the word ``metric'' means ``pseudo-Riemannian metric":  it is non-degenerate but the  requirement of positive-definiteness is relaxed; the signature can be arbitrary.
For the same reason a metric for which the curvature of the associated Levi-Civita connection obeys the Einstein condition will be called an Einstein metric (we shall not use the terminology  ``pseudo-Einstein''). 
The symbol $G$ will denote, most of the time, the Lie group $\SU(3)$,  but several discussions can often be generalized to any simple (or even semi-simple) compact Lie group.

Left-invariant metrics are fully characterized by their value at the origin of the group, and therefore by a non-degenerate bilinear symmetric form on the Lie algebra, equivalently (after having chosen an appropriate basis), by an $8\times 8$ non-degenerate symmetric matrix.
Such metrics are left-invariant by construction, but their isometry group can be larger than $\SU(3)$: as a rule it is isomorphic with  $\SU(3)\times K$ where $K$, that we call the right isometry group, is some Lie subgroup of $\SU(3)$ (more about it later).
We shall give\footnote{This was already discussed long ago, with physical applications in mind,  in \cite{RCGEF}.} a parametrization, in terms of $8\times 8$ matrices, of those left-invariant metrics for which the right isometry group is $K$, for the various possible choices of $K$.

Then, for every $K$, we shall study the Einstein condition and describe the metrics for which this condition holds.
This occurs (when it occurs) for specific values of the real parameters that enter our various parametrizations of the corresponding bilinear forms.
In some cases, \ie for some choices of $K$, our analysis is complete. Unfortunately, in some other cases we could not solve the equations in full generality, and we had to assume extra relations between the otherwise independent parameters in order to complete our study.

We shall discover, along the way, several infinite families of homogenous Einstein {\sl metrics}, but once one takes into account the action of the group of diffeomorphisms on the space of metrics giving rise, in general,  to the notion of (pseudo) Riemannian  {\sl structures}, and in particular to the notion of Einstein {structures}, these families reduce, up to scaling, to only four cases, three of which were already known: the Killing metric (which is properly Riemannian and for which $K$ is $\SU(3)$),  the so-called Jensen metric \cite{Jensen}  (which is also properly Riemannian and for which $K$ is $\SO(3)$),  a particular Lorentz metric (for which $K$ is a member -- that we call $\U(1)_I$ --  of a specific conjugacy class of $\U(1)$ subgroups), and a metric of signature $(6,2)$, that was already discovered by \cite{GibbonsLuPope}, for which $K$ is trivial. 

This paper grew up from a collection of notes whose purpose was to illustrate several concepts of (pseudo) Riemannian geometry, in particular left-invariant geometry, in a case going beyond the three-sphere $S^3$, aka $\SU(2)$, for which the study of left-invariant metrics has been thoroughly studied, long ago, in many places, hence the choice of $\SU(3)$ which is the next case in the $\SU(N)$ family. For this reason the reader will find here a section, entitled ``Miscellaneous'', with contents that have to do with left-invariant geometry, but that is not directly related with the theme of Einstein metrics: There we shall discuss for instance quadratic and cubic Casimir elements (possibly modified), Laplace-Beltrami operators, sectional curvatures, Ricci decompositions, \etc These concepts are of course standard, and we shall add nothing fundamentally new to the study of their general properties, but we shall give a number of explicit results that, we hope, will entertain the reader or trigger the interest of a few students.

As for the classification of Einstein left-invariant metrics on $\SU(3)$ with given right isometry group, or of homogeneous Einstein structures, what we can offer is unfortunately incomplete since, in several cases, we  could not solve the Einstein equation in full generality while keeping all the parameters allowed by the choice of a given right isometry group. We hope that some courageous readers will take up this study. It is of no surprise that the system of equations that one needs to solve gets more and more complex, with more and more parameters, as soon as one chooses a right isometry group that gets smaller.
For this reason, even if one can easily solve these equations by hand when the right isometry group $K$ is large enough (for instance $\SO(3)$, $\U(2)$, or $\SU(3)$ itself), the use of a computer system becomes almost compulsory for smaller groups;  for example our homogenous Einstein space with Lorentz signature involves parameters that are algebraic integers of degree $15$ with large coefficients, a manual handling of such large expressions is inefficient and prone to error.
Most calculations done here were carried out using the Mathematica software system.

The paper ends with a section devoted to possible physical applications. One of them, in particle physics, using left-invariant metrics with right isometry group $\U(2)$, was described long ago, see \cite{RCGEF}, and we shall add almost nothing to this discussion, apart from resetting the problem in a slightly more general framework. Physical applications, if any, of the  existence of pseudo-Riemannian homogenous Einstein metrics on $\SU(3)$, in particular of the one that has Lorentz signature, remain to be found.

\paragraph {Reminders.}
Every non-compact connected smooth manifold admits a Lorentz metric, and a compact connected smooth manifold admits a Lorentz metric if and only if its Euler characteristic is zero (\cite{ONeil:SRgeometry}, p. 149).
A useful corollary arises when there is a non-vanishing vector field, this implying that the Euler characteristic is zero. 
In particular, any compact parallelizable manifold, including any compact Lie group (they admit many non-vanishing vector fields !), has Euler characteristic zero.
So one a priori knows that the groups $\SU(N)$, and in particular the group $\SU(3)$, admits Lorentz metrics.
In the case of Lie groups however one does not need such general theorems to establish this property since, as already recalled, any non-degenerate symmetric bilinear form with Lorentz signature on the Lie algebra gives rise to a left-invariant Lorentz metric on the group itself, by using group translations;  these metrics are obviously homogeneous, and, in our case, have an isometry group isomorphic with $\SU(3)\times K$, where $K$ is some subgroup of $\SU(3)$. 
Usually these metrics are not Einstein metrics: this is one motivation for the search of those which are such.


\section{Metrics and curvatures}

\subsection{Killing form, inner product,  and renormalized Killing inner product}
\label{sec:KillingMetric}
Since the Lie group $\SU(3)$ is compact, the Killing form is a negative definite bilinear symmetric form on the Lie algebra $\mathfrak{su}(3)$.
Its opposite, the Killing inner product,  defines an Euclidean structure on $\mathfrak{su}(3)$ and, using left translations (for instance), a Riemannian metric on the group itself:  the Killing metric.

It is useful and standard to define the renormalized Killing form by dividing the Killing form by $2 g$, where $g$ is the dual Coxeter number.
We call $k$ the Killing inner product (so the Killing form is $-k$), and $\widehat k=k/2g$  the renormalized  Killing inner product.
For $\SU(3)$, $g=3$, therefore $\widehat k=k/6$.
Warning: Through notational abuse we also call $k$ and $\widehat k$ the corresponding bi-invariant metrics on the manifold $G=\SU(3)$.

\subsection{Several basis}
\label{sec:basis}
Let  $(e_a)$ be an arbitrary basis of the tangent space at the identity of $G$, identified with ${\mathrm{Lie}}(G)$. 
 Through notational abuse we also call  $e_a$  the corresponding left-fundamental\footnote{According to the present standard terminology they are also called {\sl left}-invariant although they commute with the right-fundamental ones (in some old references left-fundamental vector fields are called right-invariant).} vector fields obtained from the latter by letting $G$ act on itself by left multiplication.
The structure constants of the global moving frame $(e_a)$  defined by the equality $[e_a,e_b]={x_{ab}}^c \, e_c$ are identified with the structure constants of the basis  $(e_a)$ in the Lie algebra.

The dual (also called inverse) Killing metric, in the moving frame $(e_a)$ has components $k^{ab}$ and therefore reads\footnote{Here and below we use the first Einstein summation convention: an index variable
 that appears twice, once as a superscript and once as a subscript, must be summed over.} $k^{-1} =  k^{ab} \, e_a \otimes e_b$.
The Killing metric itself reads $k =  k_{ab} \,  \theta^a \otimes \theta^b$, where $(\theta^a)$ is the moving co-frame dual to $(e_a)$.
Replacing $k$ by  $\widehat k$  we have similar expressions for  the renormalized Killing metric, with $\widehat k_{ab} = k_{ab}/2g$ and  $\widehat k^{ab} = 2g\, k^{ab}$.
Since $g=N$ for $\SU(N)$, we have $\widehat k = k/6$  and   $\widehat k^{-1}= 6  k^{-1}$ for $\SU(3)$.

Let  $(X_{a})$ be a basis of ${\mathrm{Lie}}(G)$ which is orthonormal for the inner product $k$. The ordered set of vectors $\widehat X_{a} =   \sqrt{2g}  \, X_a$, in particular $\widehat X_{a} =   \sqrt{6}  \, X_a$ for $G=\SU(3)$, is then an orthonormal basis for the inner product~$\widehat k$.
We  have $k^{-1} =  \delta^{ab} \, X_a \otimes X_b$ and $\widehat k^{-1} = 2g\,  \delta^{ab} \, X_a \otimes X_b$,  where $\delta^{ab}$ is the Kronecker symbol.

Define $i\, L_a \in \, \mathfrak{su}(3)$ by the equality $i \,L_a = 2 \sqrt{3} \, X_a   =   \sqrt{2} \, \widehat X_{a} $. 
It is traditional\footnote{The hermitian matrices $\lambda_a$ are usually called Gell-Mann matrices in the physics literature.} to call $i \,\lambda_ a$ the $3\times 3$ anti-hermitian matrices that represent the $i\, L_a$  in the defining representation of $\SU(3)$.

Let $E^i_j$ be single-entry $3\times 3$ matrices. 
One obtains the $\lambda_a$  as follows~:
$ 
\lambda_ 1 = E^{1} _ 2 + E^{2} _ 1, \,
\lambda_ 2 = i (E^{2} _ 1 - E^{1} _ 2), \,
\lambda_ 3 = E^{1} _ 1 - E^{2} _ 2, \,
\lambda_ 4 = E^{3} _ 1 - E^{1} _ 3, \,
\lambda_ 5 = i (E^{3} _ 1 - E^{1} _ 3), \,
\lambda_ 6 = E^{2} _ 3 + E^{3} _ 2, \,
\lambda_ 7 = i (E^{3} _ 2 - E^{2} _ 3), \,
\lambda_ 8 = \tfrac {1} {\sqrt{3}} (E^{1} _ 1 + E^{2} _ 2 - 2 E^{3} _ 3).
$

distinct

The Lie bracket of two Lie algebra elements can be written as a matrix commutator in any chosen representation. It is standard to call $- 2 \,  {f_{ab}}^{c}$  
the (real) structure constants of the  basis $i \lambda_a$, \ie  $[i \lambda_a, i \lambda_b] =  - 2 {f_{ab}}^{c}  (i  \lambda_c) $, equivalently,
$[\lambda_a, \lambda_b] = 2 i {f_{ab}}^{c}\, \lambda_c$.  
From $Tr(\lambda_a  .  \lambda_b) = 2\, \delta_{ab}$ one obtains   $4 i \, {f_{ab}}^{c} = Tr( [ \lambda_a ,  \lambda_b ]\, . \, \lambda_c)$; using cyclicity of the trace one finds that ${f_{ab}}^{c}$ is antisymmetric in its last two indices $b,c$.

At the origin of $\SU(3)$, the left-invariant vector fields $X_a$, identified with Lie algebra elements,  are expressed as matrices  $\tfrac{i}{2\sqrt{3}} \, \lambda_a$ in the defining representation.
The structure constants of the basis $(X_a)$, which is orthonormal for $k$, are therefore equal to  $\tfrac{-1}{\sqrt{3}} {f_{ab}}^c$.
 Notice that in  the adjoint representation, the generators $i L_a$  are represented by (real antisymmetric) matrices $2 f_a$ which have matrix elements  ${2\, f_{ab}}^{c} $.

\subsection{Left-invariant pseudo-Riemannian metrics and isometry groups}
\label{sec:invariantmetrics}
\paragraph{Isometry groups.}
Isometry groups of left-invariant metrics on $\SU(3)$ are isomorphic with $\SU(3) \times K$ where $K$ is a subgroup of $\SU(3)$.
 Left-invariant metrics on a Lie group are homogeneous since the isometry group acts transitively on the manifold. 
The group $K$ can be $\SU(3)$ itself (bi-invariant metrics) or some subgroup of the maximal subgroups of the latter, which, up to conjugacy, are $\U(2)=S(\U(2) \times \U(1))$ (locally $\SU(2)\times \U(1))$ and $\SO(3)$, sometimes called the ``$\SO(3)$ principal subgroup''.  
Hence, restricting oneself to closed connected subgroups, one finds that the candidates for a (right) isometry group $K$, up to isomorphism, are the members of the following  list\footnote{The Hasse diagram of nontrivial Lie subalgebras of ${\mathrm{Lie}}(\SU(3))$, up to equivalence (conjugacy by an inner automorphism),  can be found in \cite{OFarrill}.} :  $S(\U(2) \U(1)) \sim \U(2)$, $\SU(2)$,  $\SO(3)$, $\U(1)\times \U(1)$, and $\U(1)$.

Two remarks are in order here:  1)  If, for some left-invariant metric, $K$ contains $\SU(2)$, it also contains $\U(2)$ (see, below, the paragraph called ``Parametrizations''), so $\SU(2)$ should be removed from the previous list.
2) In order to discuss left-invariant metrics up to equivalence (a notion that will be made precise later), specifying $K$ up to isomorphism is a priori not enough and, in general, one needs to specify $K$ up to conjugacy;
 however, maximal compact subgroups in a connected Lie group are all conjugate, and maximal tori are also conjugate, so only the last possibility of the above list, namely the case $K=\U(1)$, needs to be specified further.

Notice that, with the exception of $\U(1)$, specifying the type of the subgroup $K$ up to isomorphism is also enough to determine the quotient $\SU(3)/K$ as a smooth  manifold. 
Again, in the case $K=\U(1)$ one has to be more specific.
These quotients\footnote{Here we only think of these quotients as homogeneous spaces defined by the pair $(G,K)$.} are well known: we have the complex projective space  $CP^2 = \SU(3)/\U(2)$, the sphere $S^5 = \SU(3)/\SU(2)$, the irreducible symmetric space $\SU(3)/\SO(3)$ (sometimes called the Wu manifold), the flag manifold $\SU(3)/(\U(1)\times \U(1)$, and the various Aloff-Wallach spaces $\SU(3)/\U(1)$.

In order to obtain a parametrization for left-invariant metrics on $\SU(3)$ invariant under a given isometry group $\SU(3) \times K$,  the first step is to specify $K$ itself, or rather its Lie algebra, in terms of $\SU(3)$ generators;  
this is conveniently done in the defining representation, in terms of the matrices $\lambda_a$.
The second step is to impose the vanishing of the Lie derivative of an arbitrary left-invariant metric (a symmetric $8\times 8$ matrix with $8(8+1)/2 = 36$ arbitrary real constant parameters) in the direction of the generators of the Lie subalgebra of the chosen isometry subgroup $K$.
Equivalently, one can impose (or check) the equality ${r}^T \, . \, h^{-1} \,.\, r= h^{-1}$ with $r = exp(u)$ for the generators $u$ of the chosen subgroup $K$  in the adjoint representation of $\SU(3)$; one takes $u=2f_a$, for the one-parameter subgroups generated by the vectors $i\,L_a$.

Our choice\footnote{Given the isomorphy types of the right isometry groups $K$, we make here specific (but of course arbitrary) choices that define $K$ as concrete subgroups of $\SU(3)$.} 
of generators for $Lie(K)$, for the various candidates, is as follows.
For $\U(2) \sim S(\U(2) \U(1))$ (locally isomorphic with $\SU(2)  \times \U(1)$), we choose the generators $\{\lambda_1, \lambda_2, \lambda_3\}$  and  $\lambda_8$.
For $\SO(3)$ we choose the generators $\{\lambda_2, \lambda_5, \lambda_7\}$\footnote{$\SU(2)$ and $\SO(3)$ are of course locally isomorphic, but not isomorphic.}.
We identify the Cartan subgroup $\U(1) \times \U(1)$ with a fixed maximal torus of $\SU(3)$,  namely the one generated by $\lambda_3$ and $\lambda_8$. 

{\small
Let us take  $e^{i \phi} \in   \U(1)$, $k,\ell \in Z$,  and call $\U(1)_{k,\ell}$ the subgroup of $\SU(3)$ defined as the set of $3\times 3$ diagonal matrices with diagonal $(e^{i k \phi}, e^{i \ell \phi}, e^{- i (k+\ell) \phi})$.  Any one-dimensional subgroup of $\SU(3)$ is conjugate to such an $\U(1)_{k,\ell}$.
Two manifolds of the type $\SU(3)/{\U(1)}_{k,\ell}$  (Aloff-Wallach spaces) are diffeomorphic, and therefore homeomorphic, if the corresponding $\U(1)$ subgroups are conjugated in $\SU(3)$.  However  Aloff-Wallach spaces are not necessarily homeomorphic, and even when they are, they may sometimes be non diffeomorphic. This subtle problem is investigated in \cite{KreckStolzWallachSpaces}. 
Consider  $S_3$ acting on the triple $(k, l,-k -\ell)$, and identify this finite group with the Weyl group of $\SU(3)$.  Take $\sigma \in S_3$. One can show \cite{Reidegeld} that the action of the latter on  the Cartan torus changes $\U(1)_{k,\ell}$ to  $\U(1)_{\sigma(k),\sigma(\ell)}$.
It is therefore enough to assume that $k \geq \ell \geq  0$, and that $k$ and $\ell$ are co-prime (multiplying $(k,\ell)$ by an integer does not change the subgroup). One recovers the special cases  $\U(1)_Y$ for $(k,\ell)=(1,1)$, and $\U(1)_I$  for $(k,\ell)=(1,-1) \sim (1,0)$.
Another labelling possibility (that hides the roles played by $k$ and $\ell$) is to introduce a single index $\upsilon$:  up to an appropriate scaling of the generator $\tfrac{k - \ell}{2} \lambda_3+ \tfrac{k + \ell}{2} \sqrt{3} \lambda_8$, the same one-dimensional subgroup $\U(1)_{k,\ell}$ (that one can may call $\U(1)_\upsilon$) is generated by $\upsilon \, \lambda_3 + \sqrt{3} \, \lambda_8 $.}

\paragraph{Notations.}

It is traditional in physics to introduce the operators $I =  \tfrac{1}{2} \, L_3$ (isospin),  $Y =  \tfrac{1}{\sqrt{3}} \, L_8$ (hypercharge) and $Q = I+Y/2$, (electric charge).  We shall use these notations.
In the fundamental representation, where one replaces $L_a$ by $\lambda_a$, these operators $I, Y, Q$ are therefore respectively represented by the diagonal matrices $\text{diag}(1/2, -1/2, 0)$, $\text{diag}(1/3, 1/3, -2/3)$ and $\text{diag}(2/3, -1/3, -1/3)$. 
One also calls $\U(1)_I$, $\U(1)_Y$, and $\U(1)_Q$, the subgroups respectively generated by $\lambda_3$,  by $\lambda_8$ and by $Q$.
When the right isometry group is a Cartan subgroup, our above specific choice for $K$ amounts to take it equal to $\U(1)_I \times \U(1)_Y$.
Notice that the subgroups $\U(1)_Y=\{(e^{i  \phi/3}, e^{i \phi/3}, e^{i (-2) \phi/3})\}$ and $\U(1)_Q=\{(e^{i (2) \phi/3}, e^{i (-1) \phi/3}, e^{i (-1) \phi/3})\}$, with $e^{i \phi}  \in   S^1$, equivalently $e^{i \phi/3}  \in   S^1$,  respectively equal to $\U(1)_{1,1}$ and $\U(1)_{2,-1}=\U(1)_{-2,1}$, are conjugated in $\SU(3)$ by a permutation of the Weyl group $S_3$ (the triple $(-2,1,1)$ being equivalent to $(1,1,-2)$), but they are not conjugated to the subgroup $\U(1)_I$. \\
{\small Some mathematical readers could ask why physicists prefer to define $Q$ and $Y$ as before, without incorporating a multiplicative factor equal to $3$ in the definition, a choice that would indeed look more natural since irreps of $\U(1)$ are labelled by integers. The problem is that, by so doing, quarks (identified with basis vectors of the fundamental representations) would have charge $\pm 1$ and the proton (identified with a specific vector in the tensor cube of the defining representation) would have charge $+3$. However, conventionally, the latter has electric charge $+1$ (minus the charge of the electron). One could of course suggest to modify the standard terminology and redefine the notion of electric charge in such a way that the electron has electric charge $-3$, but this is not going to happen!\\
 For this reason $Y$ and $Q$ are defined as above and quarks turn out to have "fractional electric charge": $\pm 1/3$ or $\pm 2/3$.}

\paragraph{Parametrizations.}
We now  give parametrizations  for the dual metric $h^{-1}$ in the basis $(X_a)$ which is orthonormal for the Killing metric,  assuming that the isometry group of $h$ is $\SU(3)\times K$.  
Remember that in the defining representation, the vector fields $X_a$, at the origin, are represented by matrices $\tfrac{i}{2\sqrt{3}} \lambda_a$.
The reader can obtain the following results by imposing the vanishing of the Lie derivative of $h^{-1}$ with respect to the generators $\lambda_k$ of $K$, \ie\footnote{As usual, a summation over repeated indices is understood.}
\begin{equation}
\label{LieDerivativeEq}
h^{ij} (
     [\lambda_k,  \lambda_i] \otimes \lambda_j +
     \lambda_i \otimes     [\lambda_k,  \lambda_j] ) = 0
\end{equation}

For the Killing metric, $K$ is $\SU(3)$, and we have $k^{-1} =  \delta^{ab} \, X_a \otimes X_b$, in  other words, $k^{ab}$ is the unit matrix $8\times 8$.
Bi-invariant metrics are multiples of the Killing metric $k$ since $\SU(3)$ is simple, and they have the same isometry group; they read  $h^{-1} =  \alpha \,  \delta^{ab} \, X_a \otimes X_b$ where $\alpha$ is some real constant.

In the same basis $(X_a)$ the parametrization of $h^{-1}$, with components $h^{ab}$, for the other choices of $K$ specified in the list~(\ref{specificK}), reads as in the following table~(\ref{parametrizations}) (the generic parameters appearing in these expressions are arbitrary real numbers and the dots stand for $0$'s) :
\begin{equation}
\label{specificK}
\begin{array}{ccc}
& \SO(3): \{\lambda_2, \lambda_5, \lambda_7\},  \qquad \U(2):\; \{\lambda_1, \lambda_2, \lambda_3;  \lambda_8 \}, \qquad  \U(1) \times \U(1):\;  \{\lambda_3,  \lambda_8 \} \\
& \U(1)_I:\; \{\lambda_3 \},  \qquad \qquad \U(1)_Y:\;  \{\lambda_8 \}. 
\end{array}
\end{equation}

{\scriptsize
\begin{equation*}
\begin{array}{ccc}
\left(
\begin{array}{cccccccc}
 \alpha  &  .  &  .  &  .  &  .  &  .  &  .  &  .  \\
  .  & \beta  &  .  &  .  &  .  &  .  &  .  &  .  \\
  .  &  .  & \alpha  &  .  &  .  &  .  &  .  &  .  \\
  .  &  .  &  .  & \alpha  &  .  &  .  &  .  &  .  \\
  .  &  .  &  .  &  .  & \beta  &  .  &  .  &  .  \\
  .  &  .  &  .  &  .  &  .  & \alpha  &  .  &  .  \\
  .  &  .  &  .  &  .  &  .  &  .  & \beta  &  .  \\
  .  &  .  &  .  &  .  &  .  &  .  &  .  & \alpha  \\
\end{array}
\right)
&
\left(
\begin{array}{cccccccc}
 \alpha  &  .  &  .  &  .  &  .  &  .  &  .  &  .  \\
  .  & \alpha  &  .  &  .  &  .  &  .  &  .  &  .  \\
  .  &  .  & \alpha  &  .  &  .  &  .  &  .  &  .  \\
  .  &  .  &  .  & \beta  &  .  &  .  &  .  &  .  \\
  .  &  .  &  .  &  .  & \beta  &  .  &  .  &  .  \\
  .  &  .  &  .  &  .  &  .  & \beta  &  .  &  .  \\
  .  &  .  &  .  &  .  &  .  &  .  & \beta  &  .  \\
  .  &  .  &  .  &  .  &  .  &  .  &  .  & \gamma  \\
\end{array}
\right)
&
\left(
\begin{array}{cccccccc}
 \alpha  & . & . & . & . & . & . & . \\
 . & \alpha  & . & . & . & . & . & . \\
 . & . & \beta  & . & . & . & . & \zeta  \\
 . & . & . & \gamma  & . & . & . & . \\
 . & . & . & . & \gamma  & . & . & . \\
 . & . & . & . & . & \delta  & . & . \\
 . & . & . & . & . & . & \delta  & . \\
 . & . & \zeta  & . & . & . & . & \varepsilon  \\
\end{array}
\right)
\\ & & \\
K = \SO(3) 
&
K = \U(2)
&
K=\U(1) \times \U(1)
\end{array}
\end{equation*}
}
{\scriptsize
\begin{equation}
\label{parametrizations}
\begin{array}{ccc}
\left(
\begin{array}{cccccccc}
 \alpha  & . & . & . & . & . & . & . \\
 . & \alpha  & . & . & . & . & . & . \\
 . & . & \beta  & . & . & . & . & \zeta  \\
 . & . & . & \gamma  & . & \theta  & \eta  & . \\
 . & . & . & . & \gamma  & \eta  & -\theta  & . \\
 . & . & . & \theta  & \eta  & \delta  & . & . \\
 . & . & . & \eta  & -\theta  & . & \delta  & . \\
 . & . & \zeta  & . & . & . & . & \epsilon  \\
\end{array}
\right)
&
\left(
\begin{array}{cccccccc}
 \varkappa _{11} & \varkappa _{12} & \varkappa _{13} & . & . & . & . & \epsilon _1 \\
 \varkappa _{12} & \varkappa _{22} & \varkappa _{23} & . & . & . & . & \epsilon _2 \\
 \varkappa _{13} & \varkappa _{23} & \varkappa _{33} & . & . & . & . & \epsilon _3 \\
 . & . & . & \alpha  & . & \gamma  & \delta  & . \\
 . & . & . & . & \alpha  & -\delta  & \gamma  & . \\
 . & . & . & \gamma  & -\delta  & \beta  & . & . \\
 . & . & . & \delta  & \gamma  & . & \beta  & . \\
 \epsilon _1 & \epsilon _2 & \epsilon _3 & . & . & . & . & \epsilon _8 \\
\end{array}
\right)
& {}
\\ & & \\
K=\U(1)_I,  \, \text{or} \, \U(1)_{k,-k}
&
K=\U(1)_Y ,  \, \text{or} \, \U(1)_{k,k}
&
\ommit{ K=\U(1)_{k,\ell}, \, \text{with} \, k > \ell > 0  \; }
\end{array}
\end{equation}
}

The parametrization obtained for the matrices given in table~(\ref{parametrizations}), for the specific subgroups $K$ given in~(\ref{specificK}) should be understood as generic ones: obviously, for particular choices of the real parameters entering these matrices the right isometry group can be larger than $K$ (for instance by taking all the diagonal coefficients equal to $1$, and by  setting to $0$ the off-diagonal ones, one recover the Killing metric, for which $K=\SU(3)$). More generally the matrices $h^{-1}$ that obey the Killing equation~(\ref{LieDerivativeEq}) for a chosen group $K$, as specified in~(\ref{specificK}), determine left-invariant metrics for which the right isometry group is equal either to $K$ or to an over-group of $K$ that should be {\sl equal or conjugated} to one member of the list~(\ref{specificK}).

\smallskip

Remarks (proofs, using~(\ref{LieDerivativeEq}) and the commutation relations in $\mathfrak{su}(3)$, are immediate, and left to the reader): \\
$\bullet$ Invariance of a metric under $\SU(2)$ implies invariance under $\U(2)$. \\
$\bullet$ Imposing invariance under $\U(1)_{k,0}$, with $k>0$ amounts, up to conjugacy, to impose invariance under $\U(1)_{k,-k}$, and therefore gives for $h^{-1}$ the same parametrization as the one obtained when $K=\U(1)_I$.\\
$\bullet$ Invariance under any $\U(1)_{k,\ell}$, with $k > \ell > 0$, implies invariance under $\U(1)\times \U(1)$. \\
We can therefore restrict our attention to the subgroups $K$ given by the list~(\ref{specificK}).

The above parametrizations were already obtained and commented in  \cite{RCGEF} for the various choices of the subgroup $K$. In the same reference, an application to particle physics was given, namely the interpretation of the mass operator for various types of mesons in terms of the Laplacian associated to left-invariant metrics for which $\text{Lie}(K) = \mathfrak{su}(2) \oplus \mathfrak{u}(1)$. We shall come back to this discussion at the end of the present article. 

 The number of free parameters appearing in the previous expressions of $h^{-1}$ could be a priori determined by considering these metrics as coming from an $ad(K)$ invariant bilinear form at the origin of the coset space $(\SU(3) \times K)/K$, with $K$ diagonally embedded, and by reducing the isotropy action of $K$ in the tangent space at the identity ($\RR^8$) into a sum of real irreducible representations (irreps).

\paragraph{Pseudo-Riemannian structures.}
In view of using the above parametrizations to explicitly determine various curvature tensors, one wants to have as few free coefficients as possible. It is therefore useful to consider pseudo-Riemannian  {\sl structures}, rather than pseudo-Riemannian {\sl metrics}. 
The group of diffeomorphisms of a manifold acts by pullback on its space of (pseudo) Riemannian metrics.
The quotient space is, by definition, the space of Riemannian structures. 
The stabilizer of this action at a given point, \ie at a given metric, is the isometry group of this metric.
Two metrics belonging to the same orbit have conjugated stabilizers, \ie conjugated isometry groups, and each stratum (that maybe contains distinct orbits) of the obtained stratification is characterized by an isometry group, up to conjugacy.
It may also happen that distinct metrics belonging to the same orbit have the same isometry group ---we shall meet one such example in what follows.

Left-invariant metrics of signature $(p,q)$ on $\SU(3)$ can be associated with elements of $GL(8,\RR)/O(p,q)$ since they can defined by arbitrary symmetric bilinear forms of prescribed signature on the tangent space at the origin of $\SU(3)$, \ie in $\RR^8$, 
but the associated Riemannian structures are associated with points of the orbit space of the latter under the action of $Ad(\SU(3)) \subset \SO(8)$. 

Equivalence under this action generically (\ie when the right isometry group $K$ is trivial) reduces the number of free parameters from $36=(8\times 9)/2$ to $28=36-8$. For $\SU(3)\times K$ invariant metrics, with $K$ non trivial, 
one may use rotations defined by elements of $Ad(\SU(3))$ that commute with the action of $K$ to decrease the number of parameters entering the matrices of table~(\ref{parametrizations}) determined by solving equation~(\ref{LieDerivativeEq}).

For instance, if $K=\U(1)_I$, this number is reduced from $8$ to $7$:  setting  $(h^{-1})'= {r}^T . h^{-1} . r$ with $h^{-1}$ as in table~(\ref{parametrizations}), and using $r = \exp(x \  f_8)$ 
with $x=\arctan(\theta/\eta)$, one obtains a new matrix $(h^{-1})'$ of the same family that 
can be directly obtained from $(h^{-1})$ by replacing only $\eta$ by $\eta^\prime=\sqrt{\theta^2 + \eta^2}$ and the coefficient $\theta = h^{-1}_{(6,4)}=h^{-1}_{(4,6)}=-h^{-1}_{(5,7)}=-h^{-1}_{(7,5)}$, in table~(\ref{parametrizations}), by $0$.
Since $\lambda_3$ and $\lambda_8$ commute, these two matrices $(h^{-1})$ and $(h^{-1})^\prime$ define left-invariant metrics that have the same right isometry group.\\
In a similar way the number of  parameters, if $K=\U(1)_Y$, can be reduced from $14$ to $11$: the $3\times 3$ symmetric sub-matrix $\varkappa$, in the upper left corner of $h^{-1}$, can be assumed to be diagonal.

In this way one obtains respectively $1,2,3,6,7,11,28$ parameters (instead of $1,2,3,6,8,14,36$) for the choices $K= \SU(3), \SO(3), \U(2), \U(1)\times \U(1), \U(1)_I, \U(1)_Y, \{e\}$.

 A last simplification, further reducing by one the number of  parameters, is to consider metrics only up to scale, \ie metrics that differ by a constant conformal transformation (this changes the obtained curvatures by an overall multiplicative constant). 
 The number of parameters for the previous choices of $K$, once we  identify metrics that differ by equivalence {\sl and} scaling, becomes  $0,1,2,5,6,10,27$.

\paragraph{Decomposition of a bilinear symmetric form of rank $8$ on $\SU(3)$ irreps}
\label{metricdecomposition}

The group $G=\SU(3)$ acts on the vector space of $8\times 8$ symmetric matrices ---the symmetric subspace of the tensor square of the adjoint representation. This action is not irreducible and, denoting the irreps that appear in this symmetric subspace by their dimension, we have the direct sum decomposition: $36 = 1 \oplus 8 \oplus 27$, with three terms respectively associated with the  irreps of highest weights $(0,0)$, $(1,1)$, and $(2,2)$. Let us call ${{}_1h^{-1}}$, ${{}_8h^{-1}}$, ${{}_{27}h^{-1}}$, the projections of the dual metric $h^{-1}$ on these three vector subspaces.
Calling $d_{a,b,c}=\tfrac{1}{4} \ Tr(\lambda_a [\lambda_b, \lambda_c ]_{+})$ where $[\lambda_b, \lambda_c ]_{+}=\lambda_b \lambda_c + \lambda_c \lambda_b$ is the {\sl anti}-commutator, 
and $d_a$ the (symmetric) matrices with elements $(d_a)_{b,c}=d_{a,b,c}$, it is easy to show that ${{}_1h^{-1}}=\tfrac{1}{8}\, Tr(h^{-1})  \one$ and that 
${{}_{8}h^{-1}}=\tfrac{3}{5}\,  \sum_{a=1\ldots 8}\, Tr(h^{-1} d_a) \, d_a$; the last projection,  ${{}_{27}h^{-1}}$, can be obtained by difference.
Let us illustrate this decomposition by assuming that the metric $h$ belongs to the family of metrics for which the right isometry group $K$ is (at least) $\U(2)$, with the parametrization given in~(\ref{specificK}).
One obtains immediately $h^{-1} = {{}_1h^{-1}}  + {{}_8h^{-1}} + {{}_{27}h^{-1}}$, with
{\scriptsize
\begin{equation}
\label{decomposition1827}
\begin{split}
{{}_1h^{-1}}=&A \, \one \\
{{}_8h^{-1}} =& { B \sqrt 3}\, d_8 =   B \, m_8  \quad \text{with}  \quad m_8=\text {diag} \left(1,1,1,-\frac{1}{2},-\frac{1}{2},-\frac{1}{2},-\frac{1}{2},-1\right) \\
{{}_{27}h^{-1}}=& C\,   m_{27}  \quad \text{with}  \quad m_{27} =   \text{diag} \left(1,1,1,-3,-3,-3,-3,9\right) \\
\text{where}  \; &  A=\frac{1}{8} (3 \alpha +4 \beta +\gamma ), \quad B=\frac{1}{5} (3 \alpha -2 \beta   -\gamma ), \quad C=\frac{1}{40} (\alpha -4 \beta +3 \gamma )
\end{split}
\end{equation}
}

We chose to illustrate this decomposition of $h$ (actually of $h^{-1}$) in the case $K=\U(2)$, but one can do it as well 
for the other cases\footnote{The Killing metric has projection onto $\one$ only. For the $K=\SO(3)$ family (see (\ref{specificK})), the decomposition is as in (\ref{decomposition1827})
but with  $A=\frac{5 \alpha }{8}+\frac{3 \beta }{8}$, $B=0$ and $C=\frac{3 \alpha }{8}-\frac{3 \beta }{8}$,  with the same $m_8$ but with $m_{27}=\text{diag} (1, -(5/3), 1, 1, -(5/3), 1, -(5/3), 1)$.
For the Jensen sub-family (Einstein metrics, see sect.~\ref{sec: EinsteinMetrics}), one has $A=\frac{19 \alpha }{4}$, $B=0$ and $C=-\frac{15 \alpha}{4}$.}.  
One may notice, however, that such decompositions  (that do not seem to be much used)  have no reason to be compatible with the signature, or even with the non-degenerateness, of the chosen bilinear form.
Nevertheless one can consider families,  or subfamilies, of bilinear forms for which one or several of the above projections vanish. We shall come back to this possibility in the last section.

\subsection{Curvature tensors}

Expressions for curvature tensors of the Levi-Civita connection (the torsionless metric connection) defined by an invariant metric on a Lie group can be found in various places in the literature. Unfortunately these expressions are often written in a a basis (a moving frame) for which the chosen metric is orthonormal. Here we want to study various metrics while keeping the same basis. For this reason we shall give expressions of the various curvature tensors in a basis $(e_a)$ made of arbitrary left-invariant vector fields\footnote{These formulae can be found in \cite{RCAJ}.}; we call ${x_{ab}}^c$ the corresponding structure constants: $[e_a,e_b]={x_{ab}}^c\, e_c$.

The chosen metric (call it $h$) defines musical isomorphisms between a vector space and its dual; in particular, using the structure constants ${x_{ab}}^c$ and the metric coefficients $h_{ab}$ or $h^{ab}$, one can define new symbols such as $x_{abc} = {x_{ab}}^d \, h_{dc}$,  ${x^a}_{bc} =   h^{a e} \, {x_{e b}}^d \, h_{dc}$, \etc
The term  ${x^m}_ {ik} x_ {mjl}$, for instance, extracted from~(\ref{Riemann}) below, actually means  $\sum_{m^\prime, k^\prime, l^\prime} \, {x_ {m^\prime i}}^{k^\prime} {x_ {m j}}^{l^\prime}  \,  h^{m^\prime m} \, h_{k^\prime k} \, h_ {l^\prime l}$ when expressed in terms of structure constants and metric (or inverse metric) coefficients\footnote{One could write such expressions with all the indices at the same level (writing for instance $x_{mik} x_ {mjl}$) provided one uses the second Einstein summation convention, which supposes chosen a fixed metric $h$: an index variable that appears twice at the {\sl same} level, \ie twice as a superscript or twice as a subscript, should be summed over using the chosen metric or its dual.}.
Observe that the symbols $x_{abc}$ are not, in general, antisymmetric with respect to the last two indices since the metric $h$ is not assumed to be bi-invariant.

Call $R^a_{\; bcd}$  the components of the Riemann curvature tensor. The last two indices ($c$ and $d$) are the form indices, and the first two ($a$ and $b$) are the fiber indices.
Using the metric $h$, one defines $R_{a b c d}=  h_{aa^\prime}R^{a^\prime}_{\; b c d}$.
The components of the Ricci tensor are $\varrho_{bd} = R^a_{\; bad}$ and the scalar curvature is $\tau =  \varrho^d_{\; d} :=  h^{db} \varrho_{bd}$. One can also define the Einstein tensor  $\mathtt{G}= \varrho -  \frac{1}{2} \, \tau \,  h$.
One has:
 \begin{equation}
 \label{Riemann}
\begin{split}
 R_{abcd}=
\frac {1} {4}  &
(x_ {acm}\, {x_ {bd}}^m + 2 x_ {abm} \,{x_ {cd}}^m - x_ {bcm}\, {x_ {ad}}^m
  - {x_ {ab}}^m \, x_ {mcd} + {x_ {ab}}^m \, x_ {mdc} - {x_ {cd}}^m \,  x_ {mab} \\
&  + {x_ {cd}}^m \,  x_ {{mba}} +     ({x^m}_ {ac} + {x^m}_ {ca}) (x_ {mbd} + x_ {mdb}) 
  - ({x^m}_ {bc} + {x^m}_ {cb}) (x_ {mad} + x_ {mda})
\end{split}
\end{equation}

\begin{equation} 
\label{Ricci}
\varrho_{bd} =-\frac {1} {2}  x_{mbn}\,   {x_{nd}}^m - \frac {1} {2} x_{mbn} \, {x_{m d}}^n + \frac {1} {4} {x_{mnb}} \, {x^{mn}} _d -  \frac {1} {2}  (x_{mbd} +  x_{mdb}) {{x^m} _ {n}}^n
\end{equation}

\begin{equation} 
\label{ScalarCurvature}
\tau = -\frac {1} {4}  {x^{mk}}_{n} \,  {x_{mk}}^n - \frac {1} {2}  {{x}_ {m}}^{kn}\,  {x_{n k}}^m  -     {{x_m} _ {k}}^k \, {{x^m} _ {n}}^n
\end{equation}

Notice that in order to calculate the Ricci tensor for a specific left-invariant metric, one does not need to evaluate the Riemann tensor first.

In the following we shall always express the components of the curvature tensors in the basis $(X_a)$ for which the Killing metric is orthonormal: we shall take $(e_a) = (X_a)$, hence ${x_{ab}}^c = \tfrac{-1}{\sqrt{3}} {f_{ab}}^c$ in all cases.

Notice that the last term of (\ref{Ricci}) and (\ref{ScalarCurvature}), a trace, vanishes for unimodular groups, in particular for $\SU(3)$, so we can safely drop it in the practical calculations that come next.

\section{Pseudo-Riemannian homogeneous Einstein metrics on $\SU(3)$}
\label{sec: EinsteinMetrics}

As before, {\underline {\sl the}} isometry group of a left-invariant metrics $h$ on $\SU(3)$ is denoted $\SU(3) \times K$.  It is clear that any subgroup of $K$ is also {\underline {\sl a}} group of isometries of such a metric $h$.
The inverse metrics $h^{-1}$ are parametrized as in sect.~\ref{sec:invariantmetrics} but we can also incorporate an overall (constant) real scaling factor in their definition. 
The Einstein condition for the metric $h$ reads $$\varrho = \kappa \, h$$ the real number $\kappa$ being called the Einstein constant.
Equivalently, one can solve the Einstein equation $\mathtt{G} + \Lambda \, h=0$, where $\mathtt{G}$ is the Einstein tensor; $\Lambda$ is the so-called cosmological constant (although there is no cosmological interpretation in the present context!).
For Einstein metrics one has obviously $\tau=8 \kappa$ since $\text{dim}(\SU(3))=8$, moreover  $\Lambda = \kappa - \tau/2$,  therefore $\Lambda=3\kappa$.

\smallskip
{\footnotesize Remark. A pseudo-Riemannian metric on $\SU(3)$ which is left invariant and $K$-right invariant, with $K$ a Lie subgroup,  is therefore $ad_K$ invariant and passes to an $\SU(3)$-invariant pseudo-Riemannian metric on the quotient $\SU(3)/K$, but even if the metric one starts from is an Einstein metric, the metric on the homogenous space $\SU(3)/K$ has no reason to be Einstein (and in general it is not). For instance the homogeneous metrics induced on Aloff-Wallach spaces from the Killing metric on $\SU(3)$ are not Einstein (and the so-called Aloff-Wallach metrics \cite{AloffWallach} -- that are $\SU(3)$ invariant and have positive sectional curvature -- are not Einstein either), although each of these spaces admits an homogeneous Einstein metric and even a Lorentz-Einstein metric (see \cite{Wang}). 
The aim of the previous brief comment is only to stress the fact that our purpose in the present section is to study the Einstein condition for left-invariant metrics on $\SU(3)$ itself: we shall not study what happens on its quotients.
By way of contrast, however, notice that the calculations performed in this section are the same for any Lie group with Lie  algebra ${\mathrm{Lie}}(\SU(3))$, in particular for  $\SU(3)/Z_3$, which is not homotopically trivial.}

\bigskip
We now study the Einstein condition on $\SU(3)$ for the various parametrizations of the metrics for which the right isometry group is $K$, as in (\ref{specificK}), (\ref{parametrizations}), or an over-group of the latter.

\subsubsection*{\fbox{$K=\SU(3)$}}   These are the bi-invariant metrics $h = k/ \alpha $, where $k$ is the Killing metric. 
For a simple Lie group $G$, the Ricci tensor of $k$ is $\varrho = \frac{1}{4} \, k$.  It therefore defines an Einstein space with Einstein constant $\kappa = 1/4$.  Its scalar curvature is $\tau = \text{dim}(G)/4$. 

The Ricci tensor is invariant under constant scaling of the metric (a general property), 
the Einstein condition is therefore also satisfied when $k$ is scaled by $1/\alpha$, 
the Einstein constant becoming $\kappa = \alpha /4$,  with $\tau = \alpha \, \text{dim}(G)/4$; therefore $\tau = 2 \alpha$ for $G=\SU(3)$. Moreover $\Lambda=3\alpha/4$.

\subsubsection*{\fbox{$K=\SO(3)$}}  For these metrics, the Ricci tensor is diagonal, with diagonal 
$$
\left\{\frac{1}{2}-\frac{\alpha }{4 \beta },\frac{1}{24} \left(\frac{5 \alpha ^2}{\beta
   ^2}+1\right),\frac{1}{2}-\frac{\alpha }{4 \beta },\frac{1}{2}-\frac{\alpha }{4 \beta
   },\frac{1}{24} \left(\frac{5 \alpha ^2}{\beta ^2}+1\right),\frac{1}{2}-\frac{\alpha
   }{4 \beta },\frac{1}{24} \left(\frac{5 \alpha ^2}{\beta
   ^2}+1\right),\frac{1}{2}-\frac{\alpha }{4 \beta }\right\}
 $$
 The scalar curvature, for this family, is $\tau = \frac{-5 \alpha ^2+20 \alpha  \beta +\beta ^2}{8 \beta }$.
The Einstein condition gives a second degree equation, with tho real solutions, $\alpha = \beta, \kappa = \alpha/4$, the already obtained Killing metric, and another solution, the Jensen metric \cite{Jensen}:  $\beta = 11 \alpha$, with Einstein constant $\kappa = \tfrac{21}{44} \, \alpha$.
Both are properly Riemannian (signature $(8,0)$). We recover the scalar curvature $\tau = 2 \alpha$ in the first case, and find $\tau = 42\,\alpha/11$ in the second.

{\footnotesize
\begin{equation}
\label{JensenMetric}
h^{-1}=
\left(
\begin{array}{cccccccc}
 1 & . & . & . & . & . & . & . \\
 . & 11 & . & . & . & . & . & . \\
 . & . & 1 & . & . & . & . & . \\
 . & . & . & 1 & . & . & . & . \\
 . & . & . & . & 11 & . & . & . \\
 . & . & . & . & . & 1 & . & . \\
 . & . & . & . & . & . & 11 & . \\
 . & . & . & . & . & . & . & 1 \\
\end{array}
\right)
\end{equation}}

One can recover these solutions as follows, without calculating the Ricci tensor:   write $\SU(3)$ as a principal bundle with typical fiber $\SO(3)$ over the irreducible symmetric space $\SU(3)/\SO(3)$, consider a first family of metrics $h(t)$ obtained by
 dilating the Killing metric in the direction of  fibers by an arbitrary coefficient $t^2$, their scalar curvature is $\tau(h(t)) = -\frac{5 t^2}{8}+\frac{1}{8 t^2}+\frac{5}{2}$, then define a second family $\widehat{h}(t) = (1/t^2)^{3/8} \, h(t)$, the overall scaling coefficient being chosen 
  in such a way that the Riemannian volume stays constant when $t$ varies (the determinant of $h(t)$ is $(t^2)^3$). The stationary points, with respect to $t$, of the scalar curvature $\tau(\widehat{h}(t))=(t^2)^{3/8}\, \tau(h(t))$ of the metrics $\widehat{h}(t)$ are Einstein metrics \cite{Jensen}; one obtains the equation $\tfrac{d}{dt} \tau(\widehat{h}(t))) = \text{coeff} \times (t^2-1)(t^2-1/11)$, hence the solutions. 
  
 {\small
The above is a particular case of a general construction (\cite{DAtriZiller}, \cite{WangZiller}, see also \cite{RCAJ}). Assuming that both $G$ and $K$ are simple, writing $G$ as a $K$ principal bundle over $G/K$, and dilating the Killing metric of $G$ by $t^2$ in the direction of fibers, one first
obtains the following formula for the scalar curvature of the metrics $h(t)$ on $G$:   $\tau(h(t))= \tfrac{s}{2} + c\ \tfrac{k}{4} \tfrac{1}{t^2} - k(1-c)\tfrac{t^2}{4}$, where $n=\text{dim} \, G$,  $k=\text{dim}\, K$, $s=\text{dim}\,G/K$ and $c$ is the embedding coefficient of $K$ in $G$.
This result is immediately obtained by Kaluza-Klein dimensional reduction, see for instance \cite{RCAJ}, applied to this particular fibration (in this simple case one can use O'Neill formulae for Riemannian submersions with totally geodesic fibers, see \cite{Gray}, \cite{ONeill}). 
The stationary points of the scalar curvature $\tau(\widehat{h}(t))=(t^2)^{k/n}\, \tau(h(t))$ of the metrics $\widehat{h}(t)=(t^2)^{-k/n}\, h(t)$ are Einstein metrics \cite{Jensen}. 
For an irreducible symmetric pair $(G,K)$ one has $c=1-\tfrac{s}{2k}$; in that case  $\tfrac{d}{dt} \tau(\widehat{h}(t))) = \text{coeff} \times (t^2-1)(t^2-(2k-s)/(2k+s))$.
The previous results are recovered for $G=\SU(3)$, $K=\SO(3)$, using $n=8$, $k=3$ (hence $s=5$ and  $c=1/6$).}

\subsubsection*{\fbox{$K=\U(2)$}}

The Ricci tensor is diagonal, with non-zero coefficients  $\varrho_{11} =\varrho_{22}= \varrho_{33}$,  $ \varrho_{44} =\varrho_{55} =\varrho_{66} = \varrho_{77} $,  $ \varrho_{88} $, respectively given by
$$
\frac{1}{12} \left(\frac{\beta ^2}{\alpha ^2}+2\right),
\quad   
 \frac{1}{8} \left(-\frac{\beta }{\alpha }-\frac{\beta }{\gamma}+4\right),
\quad    
 \frac{\beta ^2}{4 \gamma ^2}
 $$
The Einstein condition gives only one real solution, $\alpha = \beta = \gamma$,  \ie the family of bi-invariant metrics (proportional to the Killing metric).

\subsubsection*{\fbox{$K=\U(1)\times \U(1)$}}
The non-zero components of the Ricci tensor  are $\varrho_{11} =\varrho_{22}$,  $ \varrho_{33} $, $ \varrho_{44} =\varrho_{55} =\varrho_{66} = \varrho_{77} $, $ \varrho_{88} $, and   $\varrho_{38} = \varrho_{83} $, respectively equal to 

   \begin{equation*} 
\begin{split}
 & \frac{1}{12} \left(\frac{\gamma  \delta }{\alpha ^2}+\frac{2 \alpha  \epsilon }{\zeta
   ^2-\beta  \epsilon }-\frac{\gamma }{\delta }-\frac{\delta }{\gamma }+6\right), 
   \quad
  \frac{\epsilon ^2 \left(4 \alpha ^2+\gamma ^2+\delta ^2\right)+3 \zeta ^2 \left(\gamma
   ^2+\delta ^2\right)+2 \sqrt{3} \zeta  \epsilon  \left(\delta ^2-\gamma ^2\right)}{24
   \left(\zeta ^2-\beta  \epsilon \right)^2}
   \\
& \frac{1}{24} \left(\frac{2 \alpha  \delta }{\gamma ^2}-\frac{2 \alpha }{\delta }-\frac{2
   \delta }{\alpha }+\frac{\gamma  \left(3 \beta -2 \sqrt{3} \zeta +\epsilon
   \right)}{\zeta ^2-\beta  \epsilon }+12\right), 
   \quad
\frac{\zeta ^2 \left(4 \alpha ^2+\gamma ^2+\delta ^2\right)+3 \beta ^2 \left(\gamma
   ^2+\delta ^2\right)+2 \sqrt{3} \beta  \zeta  \left(\delta ^2-\gamma ^2\right)}{24
   \left(\zeta ^2-\beta  \epsilon \right)^2}
   \\
  & \frac{-\zeta  \epsilon  \left(4 \alpha ^2+\gamma ^2+\delta ^2\right)+\beta  \gamma ^2
   \left(\sqrt{3} \epsilon -3 \zeta \right)-\beta  \delta ^2 \left(3 \zeta +\sqrt{3}
   \epsilon \right)+\sqrt{3} \zeta ^2 (\gamma -\delta ) (\gamma +\delta )}{24 \left(\zeta
   ^2-\beta  \epsilon \right)^2}
   \end{split}
 \end{equation*}  
The Einstein condition gives only one real solution, $\alpha = \beta = \gamma = \delta = \epsilon$, $\zeta = 0$, \ie the known family of bi-invariant metrics.

\subsubsection*{\fbox{$K=\U(1)_I$}}

The parametrization of a generic left-invariant metric, with $K=\U(1)_I$, involves the eight parameters ${\alpha, \beta, \gamma, \delta, \epsilon, \zeta, \eta, \theta}$ but we know that we can fix the scale $\alpha = 1$, and set the parameter $\theta$ to $0$ since different choices for $\theta$ give metrics corresponding to the same Riemannian structure (see our discussion at the end of sect.~\ref{sec:invariantmetrics}). We are left with six parameters. The Einstein condition involves one parameter more, the Einstein constant $\kappa$.
We did not solve this set of equations in full generality: we restricted our attention to the family of metrics obtained by imposing the further constraint  $\gamma = \delta$;  in that case, one of the equations implies that  $\zeta$ should vanish. 

There are five solutions (only three if one imposes $\eta \geq 0$).
The first is the Killing metric ---as expected.
The second and third solutions only differ by a sign flip in the value of the parameter $\eta$, they are properly Riemannian Einstein metrics and they are equivalent to the Jensen solution.
The last two solutions (again, they only differ by the sign of $\eta$) are Einstein metrics with a Lorentzian signature. 

The non-zero components of the Ricci tensor are $\varrho_{11} =\varrho_{22}$, $\varrho_{33}$, $\varrho_{44}=\varrho_{55}$, $\varrho_{66}=\varrho_{77}$, 
$\varrho_{7,4}=\varrho_{6,5}=\varrho_{5,6}=\varrho_{4,7}$, $\varrho_{3,8}=\varrho_{8,3}$, $\varrho_{88}$. These seven expressions are rather huge to be displayed in an article, even after setting $\theta =0$.  
As already mentioned, one can show (it is almost straightforward but cumbersome!) that the hypothesis  $\gamma = \delta$, on top of the the Einstein condition,  implies that  $\zeta$ should vanish; 
we shall therefore only display the non-zero components of the Ricci tensor and of the metric in this simpler case, which also implies that  $\varrho_{44}=\varrho_{55}$ should be equal to  $\varrho_{66}=\varrho_{77}$ and that $\varrho_{3,8}=\varrho_{8,3}$ is $0$.
Removing duplicates, we are left with five non-zero distinct components of the Ricci tensor:

 \begin{equation*} 
\begin{split}
\varrho_{11} =\varrho_{22} =& \frac{1}{12} \left(\frac{(\gamma -\eta ) (\gamma +\eta )}{\alpha ^2}-\frac{2 \alpha }{\beta }+\frac{4 \gamma ^2}{\eta ^2-\gamma^2}+8\right), \qquad
\varrho_{33} =    \frac{2 \alpha ^2+\gamma ^2+\eta ^2}{12 \beta ^2}, \\
\varrho_{44}=\varrho_{55}=\varrho_{66}=\varrho_{77}  = &\frac{1}{24} \left(\gamma  \left(\frac{8 \alpha  \eta ^2}{\left(\gamma ^2-\eta^2\right)^2}-\frac{2}{\alpha }-\frac{1}{\beta }+\frac{12 \eta ^2 \epsilon }{\left(\gamma ^2-\eta ^2\right)^2}-\frac{3}{\epsilon}\right)+12\right), \\
\varrho_{7,4}=\varrho_{6,5}=\varrho_{5,6}=\varrho_{4,7}    =&\frac{1}{24} \eta  \left(-\frac{4 \gamma ^2 (2 \alpha +3 \epsilon )}{\left(\gamma ^2-\eta ^2\right)^2}+\frac{2}{\alpha}-\frac{1}{\beta }+\frac{3}{\epsilon }\right), \\
\varrho_{88}    =&\frac{-2 \eta ^2 \left(\gamma ^2+2 \epsilon ^2\right)+\gamma ^4+\eta ^4}{4 \epsilon ^2 (\gamma -\eta )  (\gamma +\eta )}
 \end{split}
 \end{equation*}  

The dual metric $h^{-1}$ is specified by the matrix $h^{ij}$ given in (\ref{U1Jsubclassofmetrics}):
 \begin{equation} 
 \label{U1Jsubclassofmetrics}
 h^{-1}=
\left(
\begin{array}{cccccccc}
 \alpha  & . & . & . & . & . & . & . \\
 . & \alpha  & . & . & . & . & . & . \\
 . & . & \beta  & . & . & . & . & . \\
 . & . & . & \gamma  & . & . & \eta  & . \\
 . & . & . & . & \gamma  & \eta  & . & . \\
 . & . & . & . & \eta  & \gamma  & . & . \\
 . & . & . & \eta  & . & . & \gamma  & . \\
 . & . & . & . & . & . & . & \epsilon  \\
\end{array}
\right)
 \end{equation}  

 The Einstein condition reads $\varrho_{ij}= \kappa \, h_{i,j}$ where the non-zero components of the matrix $h$ are as follows:
$$
h_{1,1}=h_{2,2} = \frac{1}{\alpha },
h_{3,3}=\frac{1}{\beta },
h_{44}=h_{55}=h_{66}=h_{77}  = \frac{\gamma }{\gamma ^2-\eta ^2},
h_{7,4}=h_{6,5}=h_{5,6}=h_{4,7}    =\frac{\eta }{\eta ^2-\gamma ^2},
h_{8,8}=\frac{1}{\epsilon }
$$ 
We have five non-linear equations and five unknowns:  the five parameters $\alpha, \beta, \gamma, \epsilon, \eta$ (but one can take $\alpha=1$), and the Einstein constant $\kappa$.

\bigskip

$\bullet$ One obvious solution of the Einstein condition is obtained by setting $\eta=0$ and by taking all the other parameters equal: one recover the bi-invariant metrics. 

\bigskip

$\bullet$ Another solution, up to scale, is obtained by setting ${\alpha=1, \beta=11, \gamma=\delta=6, \epsilon=1, \zeta=0, \eta=\pm 5, \theta=0}$. See (\ref{U1JmodifiedJensenMetric}). The Einstein constant is $\kappa = 21/44$.
The metric has signature $(8,0)$. 
{\footnotesize
 \begin{equation} 
\label{U1JmodifiedJensenMetric}
h^{-1}=
\left(
\begin{array}{cccccccc}
 1 & . & . & . & . & . & . & . \\
 . & 1 & . & . & . & . & . & . \\
 . & . & 11 & . & . & . & . & . \\
 . & . & . & 6 & . & . & 5 & . \\
 . & . & . & . & 6 & 5 & . & . \\
 . & . & . & . & 5 & 6 & . & . \\
 . & . & . & 5 & . & . & 6 & . \\
 . & . & . & . & . & . & . & 1 \\
\end{array}
\right)
 \end{equation} 
}

From the metric defined by (\ref{U1JmodifiedJensenMetric}), and using the remarks at the end of sect.~\ref{sec:invariantmetrics}, one can  obtain a one-parameter family of Einstein metrics (all defining the same Einstein structure), 
for the same parameters  $\alpha_0=1$, $\beta_0=11$, $\gamma_0=6$, $\epsilon_0 =1$,  $\zeta_0= 0$, as in (\ref{U1JmodifiedJensenMetric}),  but for arbitrary values of $\theta$,  $\vert\theta\vert \leq 5$, 
 (remember that we had imposed {\it a priori} the conditions $\theta = 0$ and $\delta = \gamma$)
while also setting  $\eta= \sqrt{\eta_0^2 - \theta^2}= \sqrt{25 - \theta^2}$ in the matrix $h^{-1}$ given in table \ref{parametrizations} for the subgroup $K=\U(1)_I$. 
All these metrics have an isometry group a priori equal or conjugated to an over-group of this particular subgroup.
 
 The solution $h^{-1}$ is reminiscent of the Jensen metric: it is easy to see that the two matrices (\ref{JensenMetric}) and (\ref{U1JmodifiedJensenMetric}) are congruent; moreover, the value of $\kappa$ is the same. 
 One is therefore tempted to think that both\footnote{They are  distinct  since the symmetric bilinear forms defined by these two matrices, written in the same basis, are distinct.} metrics define the same Riemannian structure.
One could nevertheless be puzzled by the fact that the specific group $\SO(3)$ specified in the list (\ref{specificK}) does not leave invariant the metric (\ref{U1JmodifiedJensenMetric}): only the group $\U(1)_I$ of the list (\ref{specificK}), leaves it invariant
(setting $r_a = exp(f_a)$, the reader can indeed check that, for $h^{-1}$ given by \ref{JensenMetric}, the equation ${r_a}^T \, . \, h^{-1} \,.\, r_a = h^{-1}$ holds for $a=2,5,7$, whereas, for $h^{-1}$ given by (\ref{U1JmodifiedJensenMetric}), this equation holds only for $a=3$.
The right isometry group of the latter can be obtained from the same equation by taking linear combinations of the $r_a$ with arbitrary coefficients; 
one finds that this group, of type $\SO(3)$, is generated by $\{\lambda_3, \tfrac{\lambda_4+\lambda_7}{\sqrt{2}}, \tfrac{\lambda_5+\lambda_6}{\sqrt{2}}\}$. 
Although distinct from the one specified in (\ref{specificK}), it is conjugated to the latter (because $\SO(3)$ is maximal in $\SU(3)$), 
and it contains $\U(1)_I$, as it should.

\bigskip

$\bullet$ The third solution is a Lorentz metric (signature $(7,1)$).
\ommit{Notice that any inner product in the Lie algebra defines a left-invariant metric on the group, and, of course, one could start from a bilinear form with a Lorentz signature, 
but the associated Lorentzian metric would not in general give rise to curvatures obeying the Einstein condition, this is why we provide some more details on the following particular example which indeed defines an Einstein space.
 It is obtained as follows.}
 \smallskip

Let $\epsilon$ be the (unique) real root of the 15-th degree polynomial
{\scriptsize
 \begin{equation*} 
\begin{split}
&157464000 \, x^{15}+403632720 \, x^{14}-612290016 \, x^{13}-1011752856 \, x^{12}+2420977896 \, x^{11}-160395147 \, x^{10}+8214701211 \, x^9 \\
& +22205850480 \, x^8+25959494541 \, x^7+ 13520748157 \, x^6+6727192848 \, x^5+3545761995 \, x^4-307092303 \, x^3+775200861 \, x^2+1476112248 \, x+416419380
\end{split}
 \end{equation*}  }  
 
 Let $\gamma$ be the (unique) real root of the 15-th degree polynomial
 {\scriptsize
 \begin{equation*} 
\begin{split}
&1203125 \, x^{15}-5947500 \, x^{14}+27668175 \, x^{13}-91826280 \, x^{12}+247552546 \, x^{11}-578539560 \, x^{10}+1139842990 \, x^9   \\
&   -1943457696 \, x^8+2859080697 \, x^7-3567181452 \, x^6+3705721907 \, x^5-3090965208 \, x^4+1958091648 \, x^3-862410240 \, x^2+238768128 \, x-26542080
\end{split}
 \end{equation*}  }  
 
For these values of $\gamma$ and $\epsilon$, the cubic polynomial with one indeterminate $x$
{\footnotesize
 \begin{equation*} 
\begin{split}
&x^3 (3 \gamma  \epsilon +3 \gamma -12 \epsilon )+\\
&x^2  \left(\gamma ^3 (-\epsilon )+12 \gamma ^2 \epsilon -3 \gamma ^3-12 \gamma  \epsilon ^2-4 \gamma  \epsilon +12 \gamma -48 \epsilon \right)+\\
&x \left(-12 \gamma ^3 \epsilon ^2-7 \gamma ^5 \epsilon +12 \gamma ^4 \epsilon -52 \gamma ^3 \epsilon +96 \gamma ^2
   \epsilon -3 \gamma ^5-24 \gamma ^3-48 \gamma  \epsilon ^2-64 \gamma  \epsilon \right)\\
&+ (5 \gamma ^7 \epsilon -12 \gamma ^6 \epsilon +24 \gamma ^5 \epsilon -48 \gamma ^4 \epsilon +16 \gamma ^3 \epsilon +3 \gamma ^7+12 \gamma ^5)
\end{split}
 \end{equation*}  }  
has three real roots, two are negative and one is positive;  call $\eta^2$ its positive root, and call $\eta$ the positive\footnote{One can choose the negative square root as well because $\eta$ appears only in  even powers and in products $(\gamma-\eta)(\gamma+\eta)$.} square root of $\eta^2$.
Then
$$
\beta = \frac{(\gamma -\eta ) (\gamma +\eta ) \left(\gamma ^2+\eta ^2+4\right)}{-2 \left(\gamma ^2+4\right) \eta ^2+\gamma ^2 \left(\gamma ^2+4\right)+\eta ^4}
$$

$$
\kappa = \frac{\left(\gamma ^2+\eta ^2+2\right) \left(-2 \left(\gamma ^2+4\right) \eta ^2+\gamma ^2 \left(\gamma ^2+4\right)+\eta ^4\right)}{12 (\gamma -\eta )
   (\gamma +\eta ) \left(\gamma ^2+\eta ^2+4\right)}
$$

Like $\gamma$ and $\epsilon$,  the parameter $\beta$, as well as the Einstein constant $\kappa$, can be expressed as roots of  polynomials of degree 15 with integer coefficients.
\smallskip

 $\beta$ is the (unique) real root of the polynomial
{\scriptsize
 \begin{equation*} 
\begin{split}
&420959000000 \, x^{15}-1864887536000 \, x^{14}+3473091156700 \, x^{13}-3742325355930 \, x^{12}+2779023618983 \, x^{11}-1598512715722 \, x^{10}+738336195619 \, x^9 \\
& -286057154856 \, x^8+100590932418 \, x^7-32232937198 \, x^6+8922748831 \, x^5-2060272970 \, x^4+375594480 \, x^3-51335104 \, x^2+4940624 \, x-297440
\end{split}
 \end{equation*}  } 
 
 The Einstein constant  $\kappa$ is the (unique) real root of the polynomial
 {\scriptsize
 \begin{equation*} 
\begin{split}
& 75874469299200000000 \, x^{15}-194337331275110400000 \, x^{14}+301355277599416320000 \, x^{13}-332561544757530624000 \, x^{12}+282171231781966252800
   \, x^{11}\\
&   -191136024361902738240 \, x^{10}+105464748331948650048 \, x^9-47804548501070787024 \, x^8+17858543123347792128 \, x^7-5477519217851980920  \, x^6\\
& +1363429678619072700 \, x^5-269374969407033333 \, x^4+40612859877938577 \, x^3-4362120554579953 \, x^2+293255347774576 \, x-9061971967716
 \end{split}
 \end{equation*}  }

The real $\eta^2$ is the (unique) real root of the polynomial
   {\scriptsize
 \begin{equation*} 
\begin{split}
 &7237548828125 \, x^{15}+70864769531250 \, x^{14}+314655757840625 \, x^{13}+889027170133500 \, x^{12}+1845686712291930 \, x^{11}+2969194934204748 \, x^{10} \\
 &+6007481883873834
   \, x^9+14368049748482976 \, x^8+23991657392689833 \, x^7+23305737247777970 \, x^6+9939040159739877 \, x^5 \\
   &  -2269867978871308 \, x^4-3190456836365280 \, x^3+2429318649600
   \, x^2+508754442240000 \, x-6234734592000 
       \end{split}
 \end{equation*}  } 
 Both square roots of $\eta^2$ solve the equations and therefore give rise to two distinct solutions, for the same values of the other parameters. 

This Lorentzian Einstein solution is therefore obtained for a dual metric  specified by the matrix $h^{ij}$ given in~(\ref{U1Jsubclassofmetrics}), with the above values of the parameters.
Numerically,  $\eta^2 \simeq 0.0122658$, and
 \begin{equation} 
 \label{LorentzNumerical}
\{\epsilon \simeq -0.491148,\, \gamma \simeq 0.233098,\, \eta \simeq \pm 0.110751,\, \beta \simeq 1.41407,\, \zeta = 0, \, \alpha = 1\} \quad \text{and} \quad \kappa \simeq 0.121788 
 \end{equation} 
 
One can restore the $\alpha$ dependence by scaling the parameters  $\eta, \gamma, \epsilon, \beta$,  by $\alpha$.
In that case, the Einstein constant  $\kappa$ is also multiplied by $\alpha$. Remember that the Ricci tensor is invariant under a (constant) rescaling of the metric. 

\bigskip

The scalar curvature $\tau$, for the general family of metrics specified by (\ref{U1Jsubclassofmetrics}), is 
{\footnotesize
 \begin{equation} 
 \label{scalarcurvatureU1I}
\frac{8 \alpha ^2 \beta  \epsilon  \left(\gamma ^2-2 \eta ^2\right)+2 \alpha ^3 \epsilon  \left(\eta ^2-\gamma ^2\right)+\alpha  \left(6 \beta  \eta ^2
   \left(\gamma ^2-4 \gamma  \epsilon -2 \epsilon ^2\right)-\gamma ^3 (3 \beta  (\gamma -8 \epsilon )+\gamma  \epsilon )+\eta ^4 (\epsilon -3 \beta
   )\right)-2 \beta  \epsilon  \left(\gamma ^2-\eta ^2\right)^2}{12 \alpha  \beta  \epsilon  (\gamma -\eta ) (\gamma +\eta )}
 \end{equation} }
Using the previous values of parameters, one finds that $\tau$, for the Lorentz-Einstein metric, is equal to $8\kappa$, as it should. 
Numerically,  $\tau \simeq 0.974303$. Moreover, $\Lambda$, the ``cosmological'' constant, is equal to $3\kappa \simeq 0.365363$.

\smallskip

{\sl Other properties of the obtained Lorentzian Einstein metric}: 
\begin{enumerate}
\item The matrix  $h$ has seven positive eigenvalues, and one negative: its signature is Lorentzian $(7,1)$.  Using $\alpha=1$ these numerically sorted eigenvalues are {\scriptsize$(8.17347, 8.17347, 2.90825, 2.90825, 1., 1., 0.707178, -2.03605)$}.
\item The Einstein condition gives two solutions differing from one another by flipping the sign of $\eta$.
\item One can calculate the eight principal Ricci curvatures, check that they are constant (Einstein manifolds have constant Ricci curvature), all equal to $\tau/8$.
As $\tau>0$, the Ricci signature  (the signature of the Ricci quadratic form) is $(8,0)$.
\item  We already know, from the chosen parametrization, that the right isometry group of this metric is $\U(1)_I$,  the vector field $e_3$ defined by the basis vector $X_3$ being its associated Killing vector field.
\item  This Lorentzian manifold has, at every point, a cone of time-like directions. The underlying manifold, being a Lie group, is parallelizable, orientable, and it is time-orientable for this Lorentz metric.
Numerically,  $h(X_8,X_8) = - 2.03605 < 0$, the vector field $e_8$ (which is not Killing) is therefore time-like.  Notice that the Killing vector field $e_3$ is space-like. 
The integral curve of the left-invariant vector field $e_8$ is a closed time-like curve. Moreover, it is a geodesic (it is easy to show that the covariant derivative $\nabla_{e_8}\, e_8$ vanishes). The integral curve of $e_3$ is also a geodesic.
\item  One can check that this Lorentzian Einstein metric is a stationary point of the scalar curvature, when one varies the parameters while keeping the volume fixed.
This  provides another way to obtain the above solution. 
For the metrics specified by~(\ref{U1Jsubclassofmetrics}), the determinant of $h^{-1}$ is ${\sl d}=  \alpha ^2 \beta  \epsilon  \left(-2 \gamma ^2 \eta ^2+\gamma ^4+\eta ^4\right)$,
and the scalar curvature of the family of metrics ${\sl d}^{1/8} \times h$ (for which the determinant stays equal to $1$ when the parameters vary)  is $\tau / {\sl d}^{1/8}$, 
where the expression of $\tau$ in terms of the parameters $\eta, \gamma, \epsilon, \beta, \alpha$ was given in~(\ref{scalarcurvatureU1I}).
We shall only display a few curves that illustrate the stationarity property by giving plots of $\tfrac{\partial}{\partial u} \tfrac{\tau}{ {(-\sl d)}^{1/8}}$, for $u\in \{\eta, \gamma, \epsilon, \beta\}$, in a neighborhood of the found solution\footnote{The determinant being negative around the extremum that corresponds to the obtained Einstein metric (because $\epsilon < 0$), we introduce a minus sign in front of  ${\sl d}$ in ${\sl d}^{1/8}$.}.

\begin{figure}[ht]
\centering
\includegraphics[width=10pc]{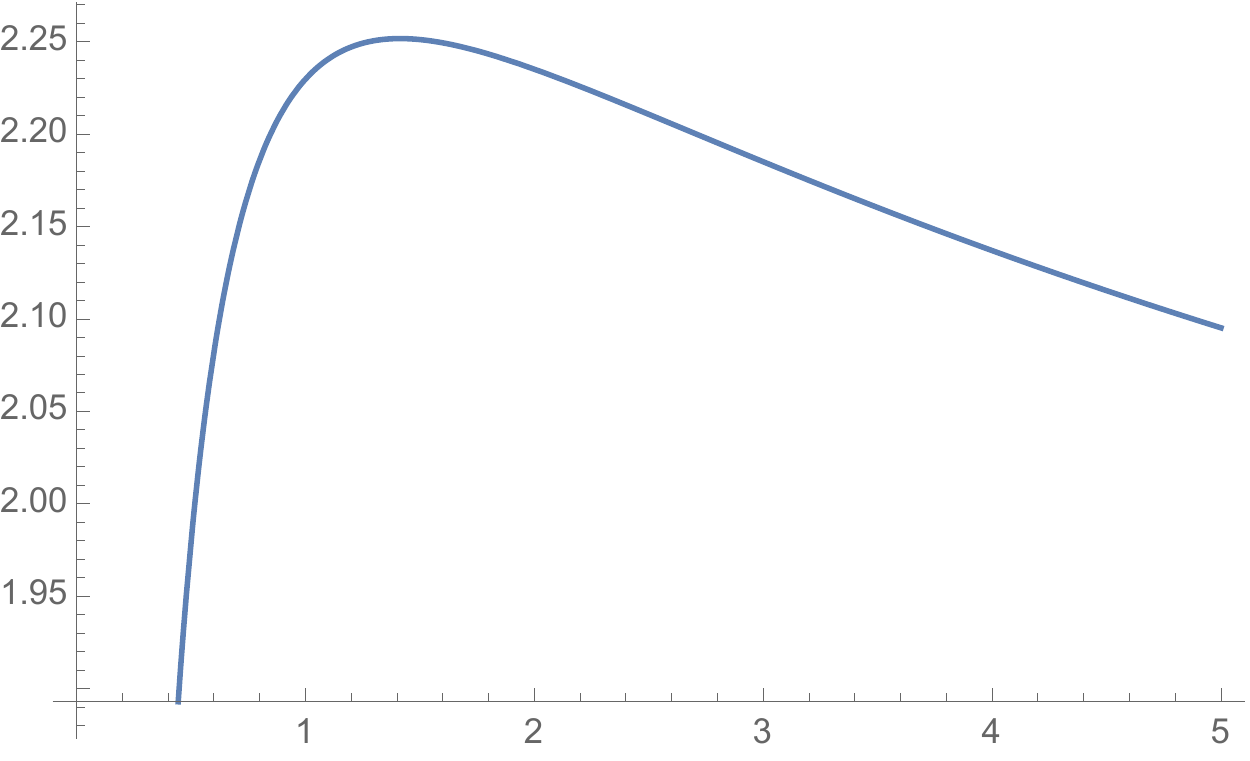}
\qquad
\includegraphics[width=10pc]{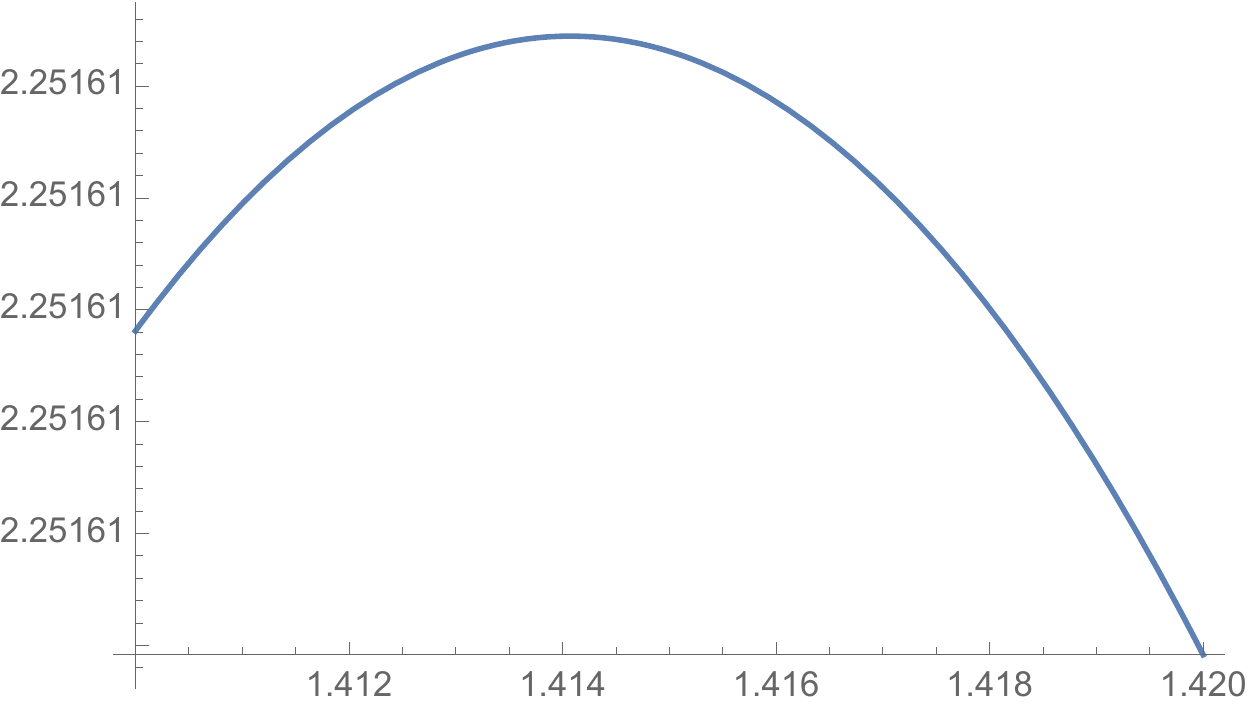}
\caption{\label{dtauoverdbeta} Derivative of $\tfrac{\tau}{ {(-\sl d)}^{1/8}}$ with respect to $\beta$, for $\beta$ in $[0,5]$ and in $[1.41, 1.42]$.}
\end{figure}

\begin{figure}[ht]
\centering
\includegraphics[width=10pc]{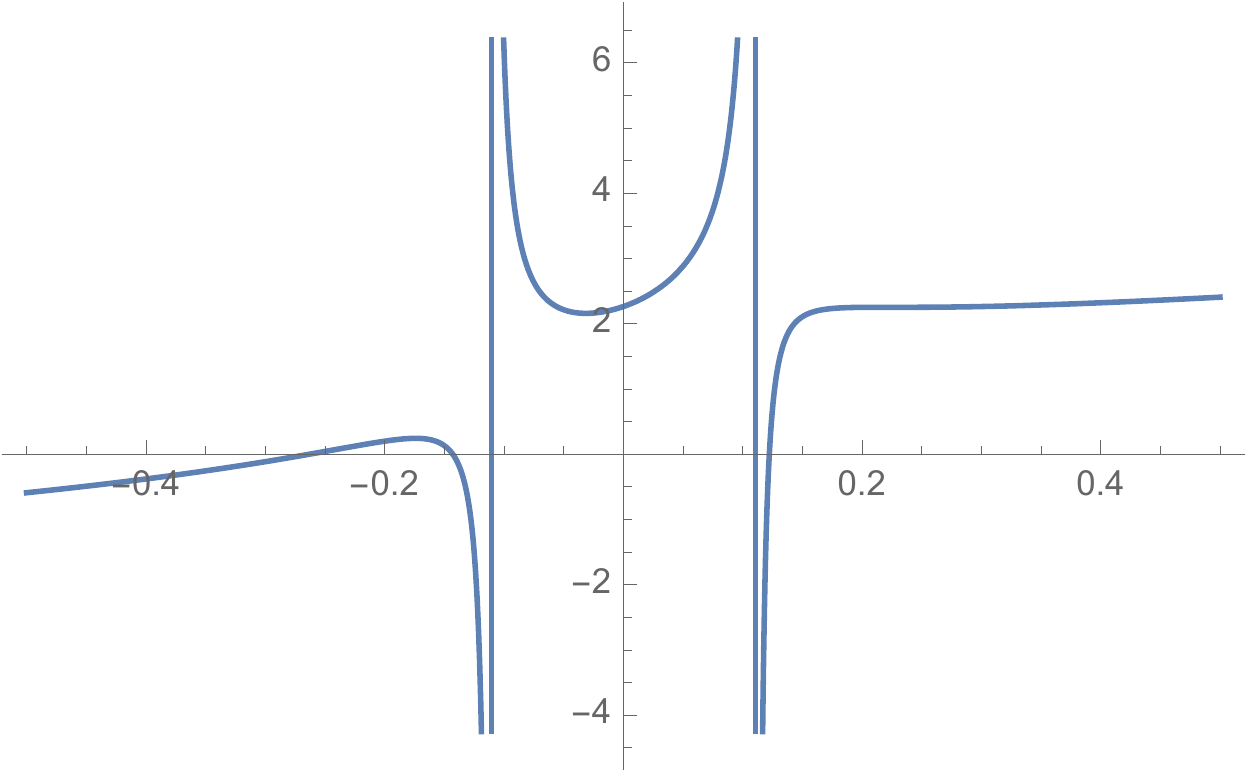}
\qquad
\includegraphics[width=10pc]{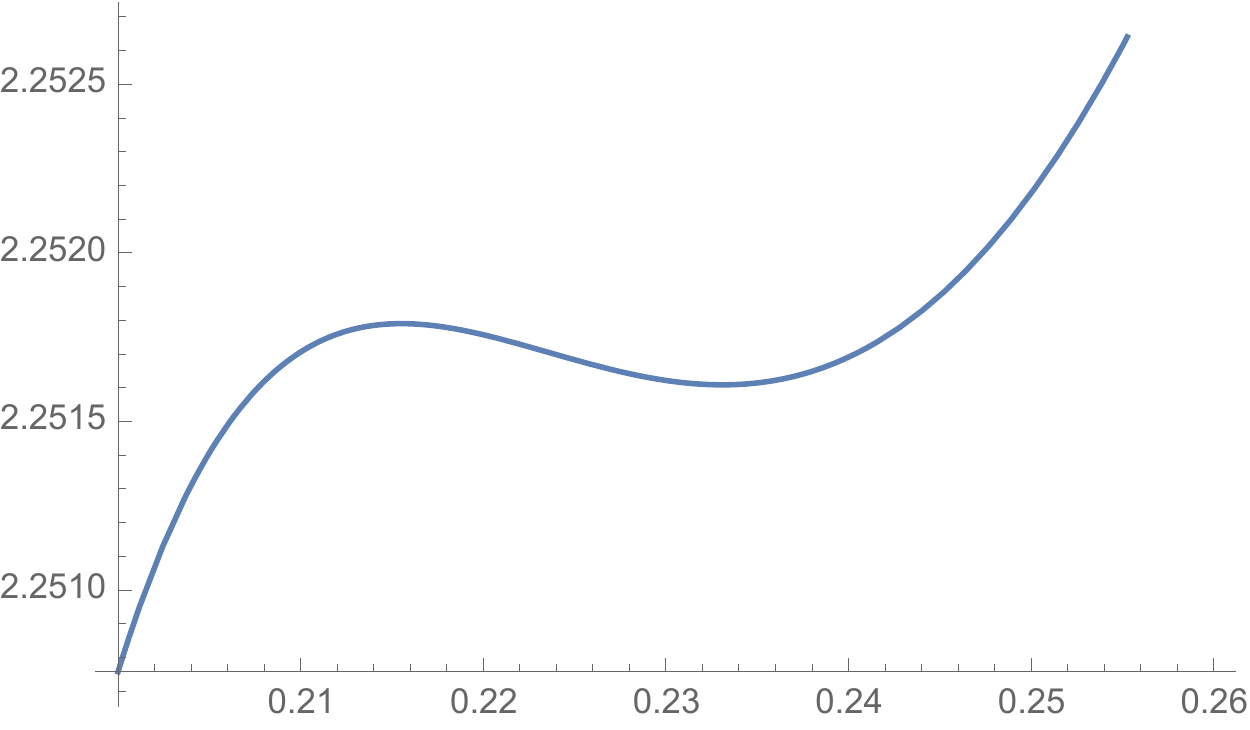}
\qquad
\includegraphics[width=10pc]{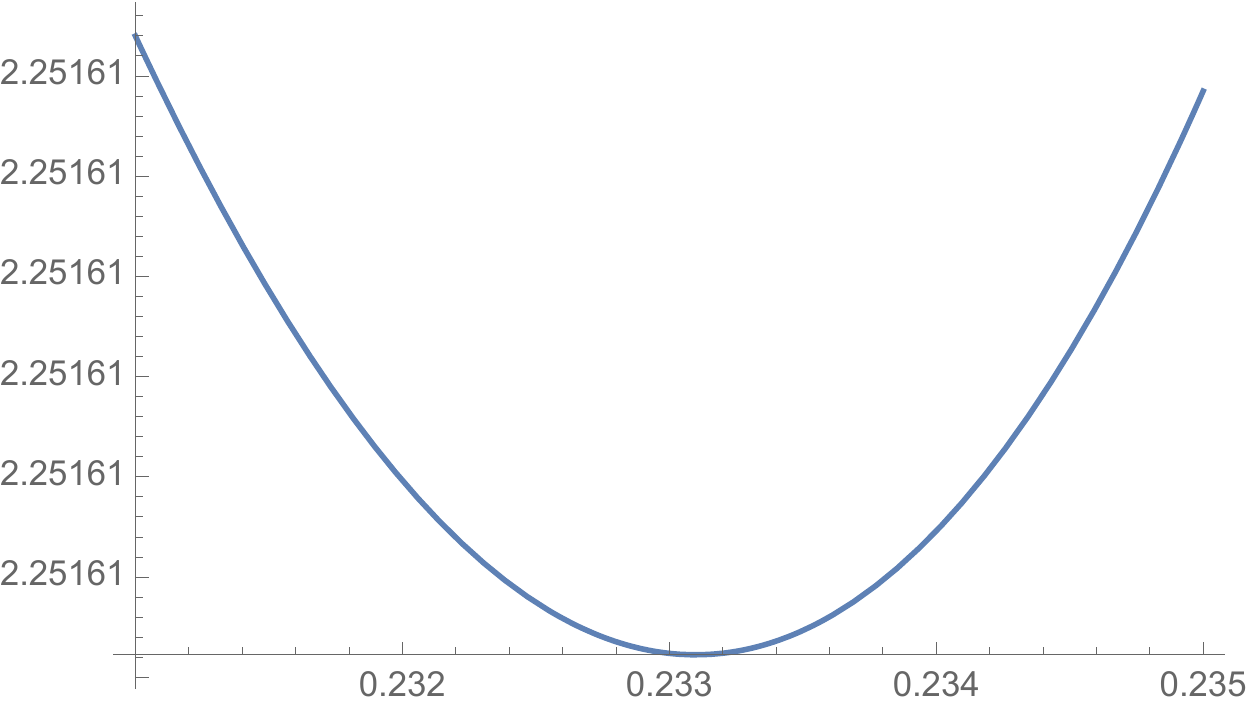}
\caption{\label{dtauoverdbeta} Derivative of $\tfrac{\tau}{ {(-\sl d)}^{1/8}}$ with respect to $\gamma$, for $\gamma$ in $[-0.5,0.5]$, $[0.2,0.26]$ and in $[0.231, 0.235]$.}
\end{figure}

\begin{figure}[ht]
\centering
\includegraphics[width=10pc]{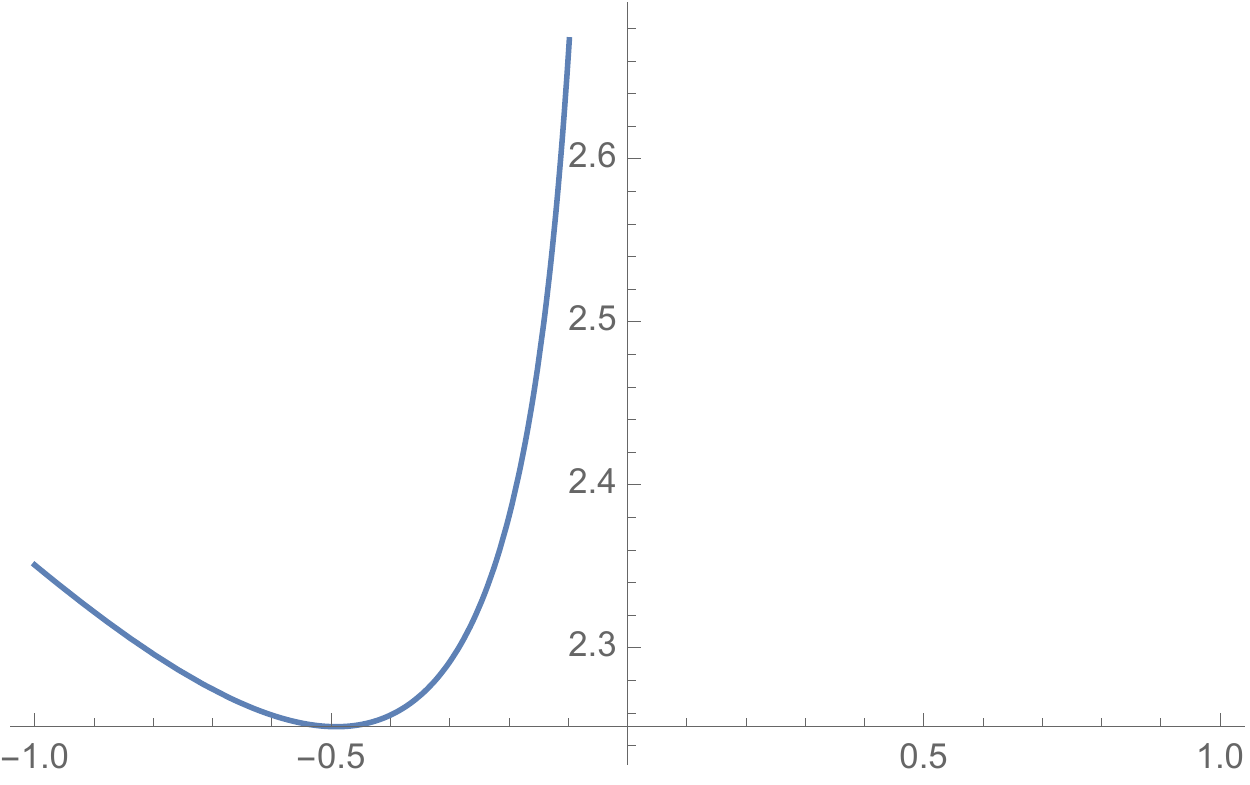}
\qquad
\includegraphics[width=10pc]{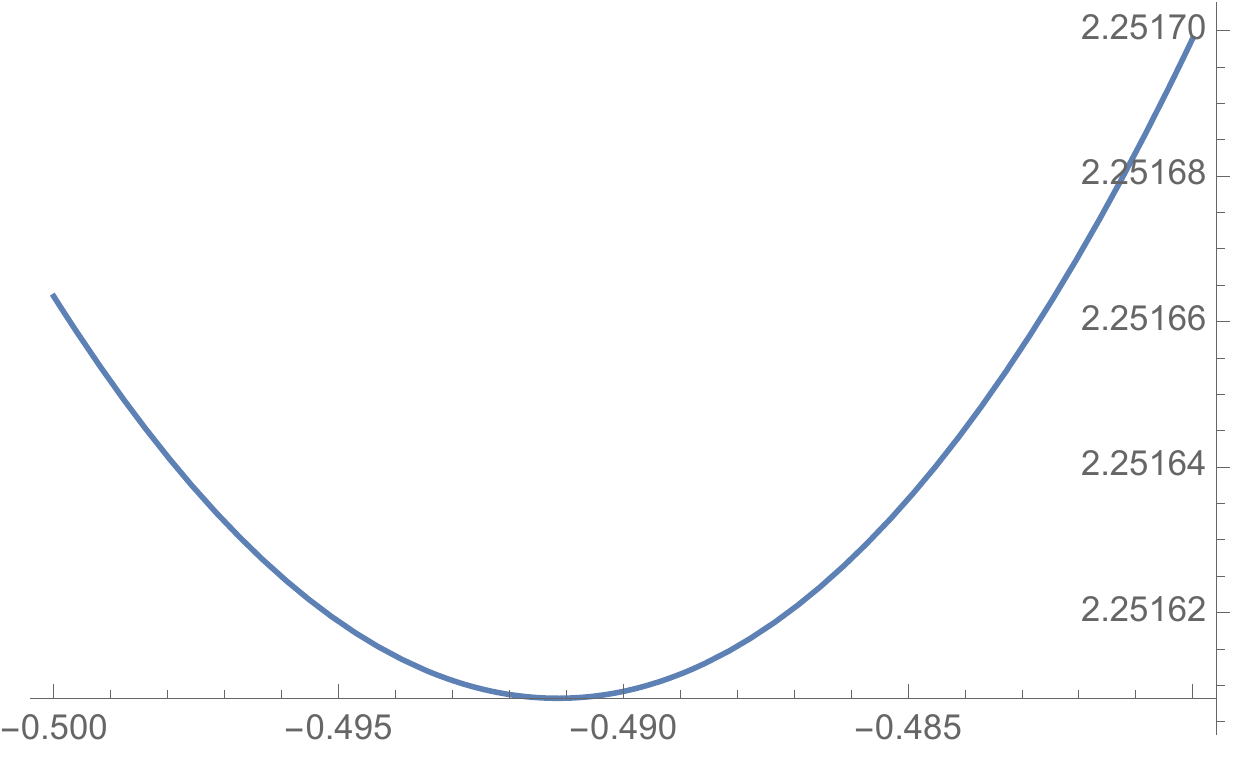}
\caption{\label{dtauoverdbeta} Derivative of $\tfrac{\tau}{ {(-\sl d)}^{1/8}}$ with respect to $\epsilon$, for $\epsilon$ in $[-1,1]$ and in $[-0.5, -0.48]$.}
\end{figure}

\begin{figure}[ht]
\centering
\includegraphics[width=10pc]{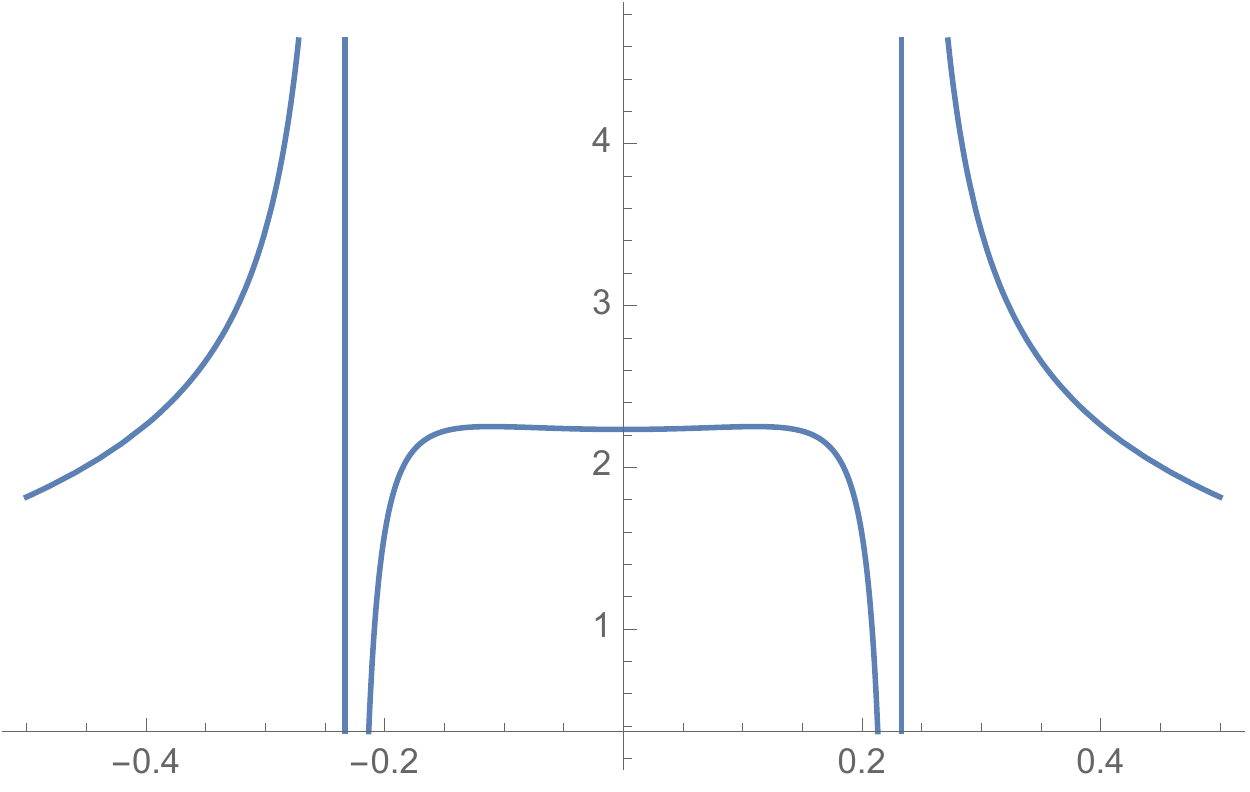}
\qquad
\includegraphics[width=10pc]{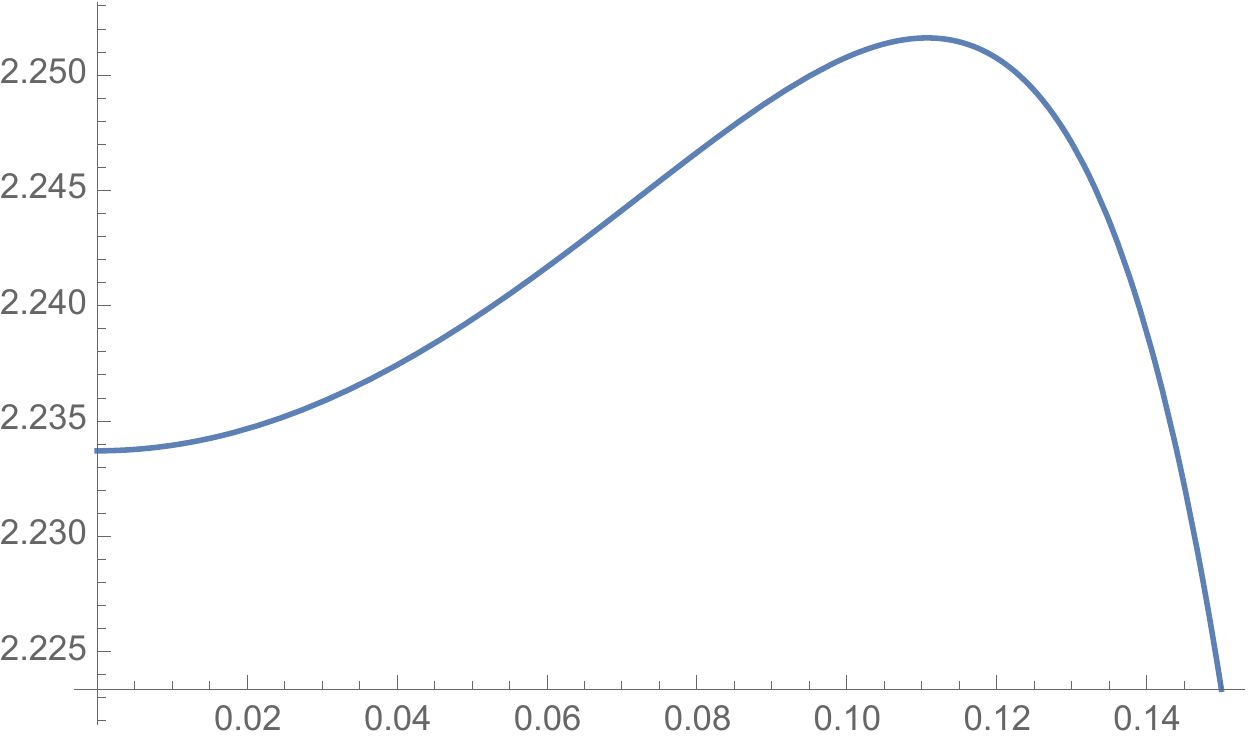}
\qquad
\includegraphics[width=10pc]{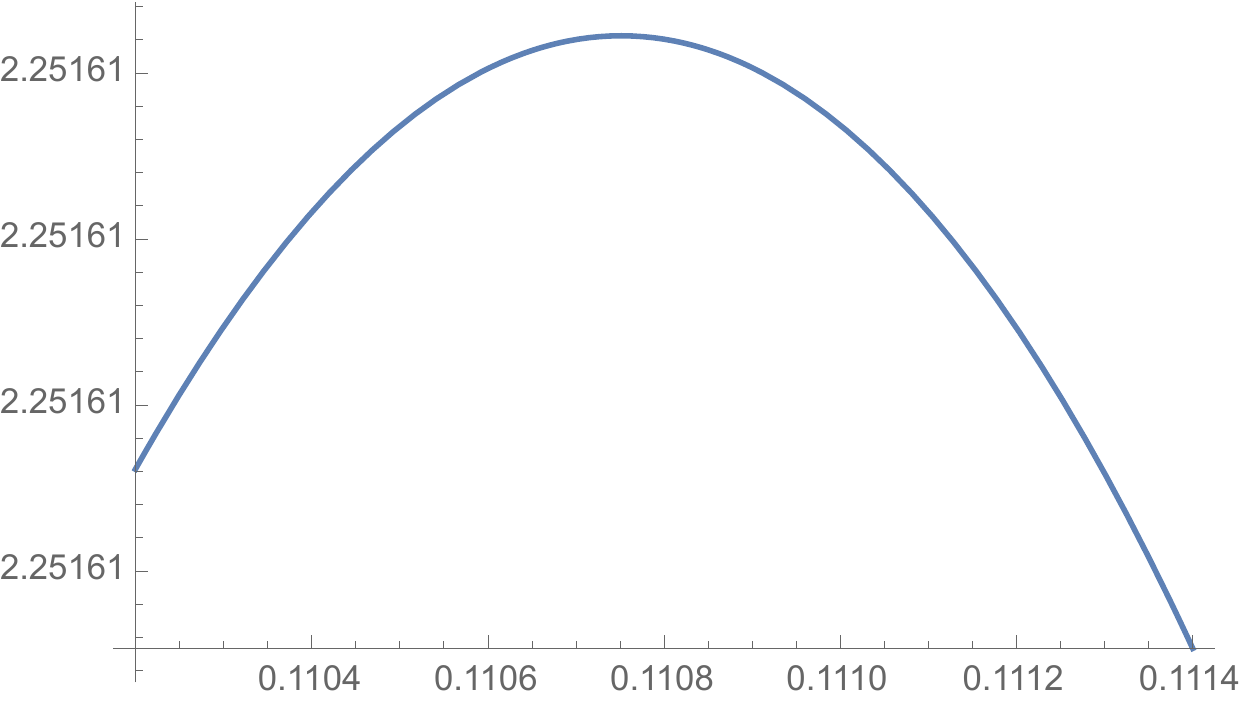}
\caption{\label{dtauoverdbeta} Derivative of $\tfrac{\tau}{ {(-\sl d)}^{1/8}}$ with respect to $\eta$, for $\eta$ in $[-0.5,0.5]$, $[0,0.15]$,  and in $[0.1102, 0.1114]$.}
\end{figure}
\item From the previous Lorentz Einstein metric, defined by parameters values that we now call $\alpha_0, \beta_0, \gamma_0$, $\epsilon_0, \eta_0,  \zeta_0$ 
 (remember that we had imposed {\it a priori} the conditions $\theta = 0$ and $\delta = \gamma$),  one obtains a one-parameter family of distinct Lorentz Einstein metrics, with the same Einstein constant, 
 for the same values $\alpha_0, \beta_0, \gamma_0, \epsilon_0, \zeta_0$, but for arbitrary values of $\theta$ (obeying $\theta^2 \leq \eta_0^2$), 
 while setting $\eta = \sqrt{\eta_0^2 - \theta^2}$ in the matrix $h^{-1}$ given in table \ref{parametrizations} for $K=\U(1)_I$ (see our discussion at the end of sect.~\ref{sec:invariantmetrics}). 
 In particular we could trade $\eta$ for $\theta$ by taking $\theta = \eta_0$, then $\eta$ vanishes.
All these metrics define the same pseudo-Riemannian {\sl structure}.
They have the same isometry group $\U(1)_I$. Notice that the calculations presented in the present subsection ($K=\U(1)_I$) do not exclude the fact that the right isometry group could be equal or conjugated to a group larger than this particular $\U(1)$, but it cannot be so, otherwise we would have already found this left-invariant Einstein Lorentzian metric in one of the previous subsections.
 \end{enumerate}

\subsubsection*{\fbox{$K=\U(1)_Y$}}
Such left-invariant metrics are parametrized by the last entry of table \ref{parametrizations}. 
Even after taking into account isometries and scaling, there are too many free parameters left  ($10$ of them) and we could not solve the Einstein condition for this family in full generality.
For this reason we looked at several subfamilies obtained by imposing conditions on the parameters, but, even then, we could not find a single example of an Einstein metric in this family, except, of course, the Killing metric, for which the right isometry group is $\SU(3)$ itself.

\subsubsection*{\fbox{$K=\{e\}$}}
Solving explicitly the system of equations coming from the Einstein condition for this family seems to be a formidable task, even after reducing the number of parameters from $36$ to $28$ by considering metrics only up to equivalence.
So we shall not have much to say in that case. 

\smallskip

We should nevertheless mention one pseudo-Riemannian Einstein metric, of signature $(6,2)$,  for which $K=\{e\}$, and that was found in \cite{GibbonsLuPope}, using other notations.
We shall describe it below.
Consider first the family of metrics defined by taking $h^{-1}$ equal to
{\footnotesize
$$ 
\left(
\begin{array}{cccccccc}
 \alpha  & . & . & . & . & . & . & . \\
 . & \beta  & . & . & . & . & . & . \\
 . & . & \gamma  & . & . & . & . & . \\
 . & . & . & \alpha  & . & . & . & . \\
 . & . & . & . & \beta  & . & . & . \\
 . & . & . & . & . & \alpha  & . & . \\
 . & . & . & . & . & . & \beta  & . \\
 . & . & . & . & . & . & . & \gamma  \\
\end{array}
\right)
$$}

For generic values of the parameters $\alpha, \beta, \gamma$, these left-invariant metrics have a trivial right-isometry group $K$ (the equation for Lie derivatives stemming from (\ref{LieDerivativeEq}) has no non-trivial solution) even though the same parameter $\beta$ occurs in positions $(2, 5, 7)$ which are those corresponding to the generators $\lambda_a$ of the $\SO(3)$ subgroup defined in (\ref{specificK}). 
For $\alpha = \beta = \gamma$ the right isometry group $K$ is $\SU(3)$, and for $\alpha = \gamma$ one recovers the cases already described in~(\ref{parametrizations}) for which $K=\SO(3)$.

The non-zero components of the Ricci tensor are $\varrho_{11} =\varrho_{44}=\varrho_{66}$, $\varrho_{22}=\varrho_{55}=\varrho_{77}$, $\varrho_{33}=\varrho_{88}$, they are respectively equal to:
$$
\left\{\frac{1}{12} \left(\frac{2 \beta  \gamma }{\alpha ^2}-\frac{\alpha }{\beta }-\frac{2 \gamma }{\beta }-\frac{2 \beta }{\gamma
   }+6\right),\frac{1}{24} \left(\frac{\alpha ^2}{\beta ^2}+\frac{4 \alpha  \gamma }{\beta ^2}-\frac{4 \alpha }{\gamma }-\frac{4 \gamma }{\alpha
   }+9\right),\frac{1}{4} \left(\frac{\alpha  \beta }{\gamma ^2}-\frac{\alpha }{\beta }-\frac{\beta }{\alpha }+2\right)\right\}
 $$
 The Einstein condition is obtained by setting the previous triple equal to  $ \{\kappa / \alpha, \kappa / \beta, \kappa / \gamma\}$.
 This system of equations has three real solutions: one first recovers the multiples of the Killing metric by taking $\alpha = \beta = \gamma$, with Einstein constant $\kappa = \alpha/4$, 
 then one recovers the multiples of the Jensen metric, for which $\alpha = \gamma$,  $\beta =11 \alpha$, and $\kappa = \tfrac{21}{44} \, \alpha$; 
  finally one obtains a third solution that we describe now.
   
  Let $P$ be a cubic polynomial with one indeterminate $x$ and real coefficients, call ${\mathfrak r}(P)$ its smallest real root. 
  Then, taking $\beta = \alpha\,{\mathfrak r}(85 x^3-29 x^2+27 x-3)$ and $\gamma= \alpha\,{\mathfrak r}(768 x^3+128 x^2+204 x+45)$ defines an Einstein metric with   
 Einstein constant $\kappa =  \alpha\, {\mathfrak r}(14400 x^3-5520 x^2+1044 x-101)$ and scalar curvature $\tau = \alpha \,  {\mathfrak r}(-808 + 1044 x - 690 x^2 + 225 x^3)$.
  Numerically $\beta/\alpha \simeq 0.121$, $\gamma/\alpha \simeq -0.213$, $\kappa/\alpha \simeq 0.196$, $\tau/\alpha \simeq 1.568$.
 This left-invariant Einstein metric has signature $(6,2)$ and its right isometry group  $K$ is trivial.
 
As already mentioned this solution was already found in \cite{GibbonsLuPope} where the authors give the matrix elements of $h$ (not of its inverse $h^{-1}$), using a different scaling, in terms of two reals $x_1, x_2$. Their values can be compared to the above ones by
writing $h = ({1}/{\beta}) \text{diag}(\{x_1, 1, x_2, x_1, 1, x_1, 1, x_2\})$; one finds $x_1= \beta/\alpha$ (given above) and $x_2= \beta/\gamma ={\mathfrak r}(768 - 128 x - 1860 x^2 + 1275 x^3)\simeq -0.570$.
Notice that $x_2 = - \tfrac{(1-x_1)(1-5x_1)}{5x_1}$. The Einstein constant for the metric $\beta \, h$ is  $\kappa / \beta \simeq 1.616$, 
and can be written\footnote{The Einstein constants given in reference \cite{GibbonsLuPope} differ from ours by two overall multiplicative factors: 
one comes from the fact that their matrix expression of $h$, compared to ours, is rescaled by $\beta$, and the other (equal to $3$)
comes from the fact that the basis vectors used by these authors to define their metrics differ from our basis vectors $(X_i)$ by a factor $\sqrt{3}$.}
 $\frac{(1-x_1) (10 x_1-1)}{20 (1-5 x_1) x_1^2}$.

\section{Miscellaneous}

\subsection{The quadratic Casimir operator}
\label{CasimirGeneral}

The quadratic Casimir element of the simple Lie group $G$ for the renormalized Killing form
 (resp. for the Killing form) is the element of the universal enveloping algebra 
 defined\footnote{We remind the reader that the Killing inner product $k$ is the opposite of the Killing form, hence the minus sign in front of the expressions defining $\Omega_2^{k}$ and ${\widehat \Omega_2}$, since $(X_a)$ is an orthonormal basis for $k$. 
 See sect.~\ref{sec:basis}.} 
 by ${\widehat \Omega_2}= - \sum_a \widehat X_{a} .  \widehat X_{a}$ (resp. ${\Omega_2}= - \sum_a X_a . X_a$). 
Casimir elements can be evaluated in any representation, and, in an irreducible representation,  ${\widehat \Omega_2}$ (resp. ${\Omega_2}$) is a multiple of the identity matrix, with eigenvalue ${\widehat C_2}$ (resp. ${C_2}$).
The definition of Casimir operators involves the {\sl inverse\/} Killing inner product, so, using  $\widehat k=k/2g$,  one obtains the relation\footnote{We also remind the reader that $g$ is the dual Coxeter number, which is equal to $N$ for $\SU(N)$.}: 
\begin{equation}
\label{Casimir2krenorm}
{C_2} = {\widehat C_2} / 2g.
\end{equation}

Explicitly, for an irreducible representation of highest weight $\w$,
one obtains 
\begin{equation}
\label{Casimir2K}
{\widehat C_2} =   \langle \w+\rho, \w +\rho  \rangle  -   \langle \rho, \rho \rangle =  \langle \w, \w+ 2 \rho \rangle
\end{equation}
where $\rho$ is the Weyl vector and $\langle . , . \rangle$ is the Cartan inner product in the space of roots, normalized in such a way that the length square of long roots is equal to $2$. 
One has also:
\begin{equation}
\label{Casimir2k}
{C_2} =  \sum_{\alpha}  \langle \w+\rho, \alpha \rangle^2  -   \langle \rho, \alpha \rangle^2
\end{equation}
where $\alpha$ runs over the set of all roots (use the identity $\sum_{\alpha} \vert \alpha \rangle \langle \alpha \vert =2g$ to relate (\ref{Casimir2K}) and (\ref{Casimir2k}) as in (\ref{Casimir2krenorm})).

For $\SU(N)$  in the defining representation one obtains ${\widehat C_2}=  (N^2-1)/N$ and ${C_2} = (N^2-1)/2$. In the adjoint representation  one obtains ${\widehat C_2}= 2N$ and ${C_2}= 1$.

In the case of $\SU(3)$, one can use for instance (\ref{Casimir2K}) to show that, for an irreducible representation of highest weight $\w$ with (Dynkin) components $(o_1, o_2)$ in the basis of fundamental weights,
\begin{equation}
\label{Casimir2KSU3}
{\widehat C_2} =  \frac{2}{3} (o_1^2 + o_1 o_2 + o_2^2) + 2 (o_1 +  o_2)
\end{equation}
Equivalently, one can  evaluate ${\widehat \Omega_2}= - \sum_a  \tfrac{i L_a}{\sqrt{2}} . \tfrac{i L_a}{\sqrt{2}}$  and  ${\Omega_2} = -  \sum_a  \tfrac{i L_a}{2\sqrt{3}} . \tfrac{i L_a}{2\sqrt{3}}$ in the chosen representations.
The above general relations, in the case of $\SU(3)$, give:   ${\widehat C_2}=6$, ${C_2}=1$ in the adjoint representation, and ${\widehat C_2}=8/3$, ${C_2}=4/9$  in  the defining representation  (these values can also be directly calculated by representing the $i L_a$ generators  by matrices $2 f_a$ in the former case and by matrices $i \lambda_a$ in the latter).
\smallskip

{\small For the group $\SU(2)$, and for an irreducible representation of highest weight $2 j$ (where  the ``spin'' variable $j$ is an integer or a half-integer), of  dimension $2j+1$, the value  $j(j+1)$ presented in the majority of quantum physics textbooks as eigenvalue of ``the Casimir operator'' corresponds to a Casimir element neither associated with the Killing form on $\SU(2)$ (${C_2} = j(j+1)/2$) nor with the renormalized Killing form (${\widehat C_2} = 2j(j+1)$). 
Details: the unique long root, which is also the highest weight $\sigma = 2$ of the vector representation (of dimension 3), obeys $<2, 2> = 2$,  so  $<1, 1> = 1/2$,
and (\ref{Casimir2K}), using $\rho =  1$,  indeed gives  ${\widehat C_2} \, = \, <2j + 1, 2j+1> - < 1,1> = ((2j+1)^2-1) <1, 1>=4 j  (j+1)  <1, 1>= 2j(j+1)$.  In order to obtain $j(j+1)$ one has to use another rescaled Killing form,  namely $k/2 = 2 {\widehat k}$, in which case the associated Casimir can still formally be given by the rhs of \ref{Casimir2K}, provided one normalizes the Cartan inner product in such a way that the length square of long root is equal to $1$, a choice that is also often made in the same quantum physics textbooks (but remember that for us this length square is equal to $2$).

\paragraph{Dynkin index.}
 In  an arbitrary basis $(e_a)$, we  have $Tr(\w(e_a) \w(e_b)) =  - 2  \iota_\w  \; {\widehat k}_{ab} =  - \iota_{\w}/g  \, k_{ab}$.
 Here $\iota_{\w}$ denotes the Dynkin  index\footnote{Some authors incorporate a pre-factor $2$ in the definition of the Dynkin index.}  of  the representation $\w$ of the Lie group $G$.
 If $\w$ is the defining representation of $\SU(N)$, one has $\iota_\w=1/2$.  If $\w$ is the adjoint representation of $G$, one has $\iota_\w=g$; in particular, for $G=\SU(N)$,  $\iota_\w=N$.
More generally, one has the relation:  ${\widehat C_2}   =   2 \,  \iota_{\w} \times {\text{dim}(Lie(G))}/{\text{dim}(\w)}$.

\subsection{Restriction to subgroups: branching}
\label{RestrictionSubgroup}

We consider the Lie algebra embedding ${\mathrm{Lie}}(\U(2)) \subset {\mathrm{Lie}}(\SU(3))$ \ie  $\mathfrak{su}(2)\oplus \mathfrak{u}(1)  \subset \mathfrak{su}(3)$, and we take $\mathfrak{u}(1)$ as the Lie algebra of the subgroup called $\U(1)_Y$ in previous sections.
This is a Levi type subalgebra : the set of simple roots of the semi-simple component of the subalgebra can be chosen as a subset of the set of simple roots of the given Lie algebra.
Call $\alpha_1, \alpha_2$ the simple roots of $\mathfrak{su}(3)$ and $\omega_1, \omega_2$ its fundamental weights.
We take $v=\alpha_1$ as the simple root of $\mathfrak{su}(2)$ (the ``$v$'' stands for ``vector'' since the $\mathfrak{su}(2)$ irrep of highest weight $v$ is the vector representation) and $t$ the fundamental $\mathfrak{u}(1)$ weight. 
The $\U(1)_Y$ generator is $3\, Y = \sqrt{3}\, L_8$ and reads $\sqrt{3}\, \lambda_8= \text{diag} (1,1,-2)$ in the defining representation; its eigenvalues are integers, as they should.
Notice that  $k(\sqrt{3}\, i L_8,\sqrt{3}\, i L_8)= 3 \times12$, so $\widehat k(\sqrt{3}\, i L_8,\sqrt{3}\, i L_8)= 3 \times12/6 = 6$ and $\widehat k^{-1}(t,t)=1/6$.
Notice also that $v=2\sigma$, where $\sigma$ denotes the $\mathfrak{su}(2)$ fundamental weight\footnote{The component along $\sigma$ of each weight of an irrep of $\mathfrak{su}(2)$ is equal to twice the ``(iso-)spin''. 
For instance those of the spinorial irrep (highest weight $\sigma$) are twice $\pm 1/2$, those of the vectorial irrep (highest weight $v$) are twice $(1,0,-1)$.}.

The simple root $\alpha_2$  of $\mathfrak{su}(3)$  is a priori a linear combination of $v$ and $t$: we have $\alpha_2= a \, v + b \, t$.
We determine $a$ and $b$ from the inner products of roots and weights calculated using the Cartan matrix or its inverse. 
As usual, all roots have length $2$ both for $\mathfrak{su}(3)$ and for $\mathfrak{su}(2)$ (we have only long roots here), so $\langle \,  \alpha_1, \alpha_1  \rangle = \langle \,  v, v  \rangle = 2$.
From the Cartan matrix of $\mathfrak{su}(3)$, namely $\begin{psmallmatrix} 2 & -1\\-1 & 2\end{psmallmatrix}$,
we get  $\langle \,  \alpha_1, \alpha_2  \rangle = -1$, moreover $\mathfrak{su}(2)$ and $\mathfrak{u}(1)$  are  orthogonal subspaces for $\langle \, , \,  \rangle$, so $\langle \, v, \,t  \rangle=0$, therefore 
$a\langle \, v, \,v \rangle=-1$, and we obtain $a = -1/2$. 
We have also  $\langle \,  \alpha_2, \alpha_2  \rangle = 2$, therefore  $a^2  \langle \,  v, v  \rangle  +     b^2  \langle \,  t, t  \rangle = 2$.
Using  $\langle \,  t, t  \rangle=1/6$ one gets $b=3$. Therefore $\alpha_1=v$, and $\alpha_2= - v/2 + 3 t$. 

The restriction matrix defining the embedding in terms of fundamental weights (which also gives the $\U(2)$ weight components $2I$ and $3Y$ from the Dynkin components $(o_1, o_2)$ of the highest weight $\w$ of any irreducible $\SU(3)$ representation) reads:
\begin{equation}
\label{restrictionmatrix}
\left(
\begin{array}{c}
 \omega_1 \\
 \omega_2 \\
\end{array}
\right)
=
\left(
\begin{array}{cc}
 1 & 1 \\
 0 & 2 \\
\end{array}
\right)
\left(
\begin{array}{c}
 \sigma \\
 t \\
\end{array}
\right)
\qquad
(2I, 3Y) = \left(o_1, o_2\right) \, \left(
\begin{array}{cc}
 1 & 1 \\
 0 & 2 \\
\end{array}
\right)
\end{equation}

Examples.\\ 
Consider the basic (fundamental) irrep of $\SU(3)$ with highest weight $\w=(1,0)$, of dimension $3$. Using the restriction matrix (eq~\ref{restrictionmatrix}) on the weight system of $\w$, namely $\{(1, 0), (-1, 1), (0, -1)\}$ we obtain the weights  appearing in the branching from 
 $\mathfrak{su}(3)$ to $\mathfrak{su}(2)\oplus \mathfrak{u}(1)$, namely $\{(1, 1), (-1, 1),  (0, -2)\}$; the associated decomposition of irreps, in terms of highest weights, reads $(1,0) \rightarrow  (1,1) \oplus (0,-2)$ 
 where, on the right hand side, the first member ($2I$) of each pair is the component along $\sigma$ of the $\SU(2)$ highest weight and where the second member ($3Y$) is the component of the $\U(1)$ weight along $t$.
Equivalently, in terms of dimensions\footnote{Remember that an $\SU(2)$ irrep with highest weight $2I$ (\ie spin $I$) has dimension $2I+1$, and that an $\SU(3)$ irrep with highest weight components $(o_1,o_2)$ has dimension  $(o_1+1)(o_2+1)(o_1+o_2+2)/2$.}:  
$[3]  \rightarrow [2]_1 \oplus [1]_{-2}$  where the subindex of $[2I+1]_{3Y}$ refers to the component of the $\U(1)$ weight.
 Conservation  of the $\U(1)$ (hyper) charge reads $2\times(1)+1\times(-2)=0$. \\
 For the adjoint representation (highest weight $\w=(1,1)$, of dimension $8$), the branching rule can be obtained in the same way and reads, when written in terms of dimensions (no confusion can arise in this case):  $[8]  \rightarrow [3]_0 \oplus [2]_3 \oplus [2]_{-3} \oplus [1]_0$. \\
 Let us conclude this section with a slightly more involved example: we consider the $\SU(3)$ representation of highest weight $\w=(2,1)$, which is of dimension $[15]$. Using the restriction matrix 
 on the weight system\footnote{$\{ (2, 1), (3, -1), (0, 2), (1, 0), (1, 0), (-2, 3), (2, -2), (-1, 1), (-1, 1), (0, -1), (0, -1), (-3, 2), (1, -3), (-2, 0), (-1, -2) \}$.}
of this highest weight of $\SU(3)$, we obtain the weights\footnote{$\{ (2, 4), (3, 1), (0, 4), (1, 1), (1, 1), (-2, 4), (2, -2), (-1, 1), (-1, 1), (0, -2), (0, -2), (-3, 1), (1, -5), (-2, -2), (-1, -5) \}$.} appearing in the branching to $\U(2)$.
The associated decomposition reads $(2,1)\rightarrow (3,1) \oplus (2,4) \oplus (2,-2) \oplus (1,1) \oplus (1,-5) \oplus (0,-2)$ where, again, on the right hand side, the first member ($2I$) of each pair is the component along $\sigma$ of the $\SU(2)$ highest weight and where the second member ($3Y$) is the component along~$t$ of the $\U(1)$ weight.
In terms of dimensions, this rhs reads $[4]_1 \oplus [3]_4\oplus [3]_{-2}\oplus [2]_1 \oplus [2]_{-5} \oplus [1]_{-2}$ and we can check the conservation  of the $\U(1)$ (hyper) charge: $4\times(1)+3\times(4) + 3\times (-2)+ 2\times(1) + 2\times(-5) + 1 \times (-2)=0$. 

\subsection{Laplacian}
Let $h$ be a Riemannian or pseudo-Riemannian metric  on the Lie group $G$.  
Assuming that $h$ is left-invariant (hence homogeneous), we can write  $h = h_{a b} \, \theta^a \otimes \theta^b$ where $h_{ab}$ are constants (real numbers) and $(\theta^a)$ is the global moving co-frame dual to the arbitrary moving frame $(e_a)$ defined from an arbitrary basis, also called $(e_a)$, in the Lie algebra  of $G$ identified with the tangent space to $G$ at the identity. The dual (\ie inverse) metric reads $h^{-1} = h^{a b} \, e_a \otimes e_b$. 
With the usual convention, the rough metric Laplacian (or Laplace-Beltrami operator)  on functions on the manifold has negative  spectrum ---so it is the {\sl opposite} of  the De Rham Laplacian on $0$-forms--- and can be written as the second-order differential operator $\Delta =    h^{a b} \, e_a \circ e_b$ where the vector fields $e_a$ act on functions on $G$. More generally, when studying the action of the Laplacian on sections of vector bundles over $G$,  the $e_a$ would act as a Lie derivative of sections in the direction~$a$.

\paragraph {Laplacian of bi-invariant metrics.}

We  call $\Delta_0$ the Laplacian associated with the Killing metric $k$; its eigenstates are labelled by irreducible representations $\w$ of $G$, the eigenvalues of $-\Delta_0$ are equal to the Casimir eigenvalues ${C_2}$ (see \ref{Casimir2k}) evaluated in the representation $\w$, and the degeneracy is $dim(\w)^2$,  see \cite{Beers} and \cite{Fegan}.

\paragraph {Laplacian of left-invariant metrics.}

Let $h$ be an arbitrary left-invariant metric on the Lie group $G$,  the spectrum of the corresponding Laplacian  $\Delta$ is discussed in a number of places (see for instance \cite{Lauret_su2}, \cite{Lauret}, and references therein).
Using the Peter-Weyl theorem together with left-invariance of the metric one can replace a difficult problem of analysis on manifolds by a simpler algebraic problem:  as in the bi-invariant case, the eigenvalues of the Laplace operator can be obtained, up to sign, as eigenvalues of some appropriate metric-dependent modified Casimir operator (consider for instance the expression  (\ref{LaplacianU2}) below) evaluated in irreducible representations of $G$.  One should be careful with this terminology because the associated modified Casimir elements (that can be defined in the enveloping algebra of $\text{Lie}(G)$) are not, in general, central.

For left-invariant metrics $h$ with isometry group $G\times K$ and more generally for naturally reductive metrics on Lie groups one can certainly write general results but here we only want to focus on the $K$ dependence of the eigenvalues in a few specific cases, and we shall be happy with some elementary calculations. 

So we return to the case $G=\SU(3)$, call $X_a$ the vectors of an orthonormal basis for the Killing metric, $h^{ab}$ the covariant components of some chosen left-invariant metric $h$ in the same basis, and 
$C(\w)$ the list of eigenvalues of the operator $\Delta =  h^{ab} \, X_a.X_b$ evaluated in some chosen non trivial representation $\w$ of $G$. The degeneracy of each eigenvalue is at least $dim(\w)$. With the notations of sect.~\ref{sec:basis}, namely setting $X_a=\frac{i}{2\sqrt 3}L_a$, we can also write
$\Delta = -\tfrac{1}{12} h^{ab} \, L_a.L_b$.  
Taking for instance $h=k$, the Killing metric, we have
\begin{equation}
\begin{split}
\Delta_0=&\left(X_1.X_1+X_2.X_2+X_3.X_3\right) + \left(X_4.X_4+X_5.X_5+X_6.X_6+X_7.X_7\right) + X_8.X_8 \\
{} =& \frac{-1}{12}\left(\left(L_1.L_1+L_2.L_2+L_3.L_3\right) + \left(L_4.L_4+L_5.L_5+L_6.L_6+L_7.L_7\right) + L_8.L_8\right) \\
\end{split}   
\end{equation}
and replacing $L_a$ by $\lambda_a$ (in the defining representation), or by $-2 i f_a$ (in the adjoint), one recovers the known Casimir eigenvalues.

\bigskip
Let us now choose a metric $h$ for which $K=\U(2)$, with parameters $\alpha, \beta, \gamma$ as in (\ref{parametrizations}). The Laplacian reads as follows, and we may introduce the notation ${\Omega_2^{\U(2)}}$ to denote the ``modified Casimir operator'' defined as $-\Delta$.
\begin{equation}
\label{LaplacianU2}
\begin{split}
\Delta =& \frac{-1}{12}\left(\alpha \left(L_1.L_1+L_2.L_2+L_3.L_3\right) +\beta  \left(L_4.L_4+L_5.L_5+L_6.L_6+L_7.L_7\right) +\gamma L_8.L_8\right)
\end{split}   
\end{equation}

In the fundamental representation $\w=(1,0)$ of $\SU(3)$, $-\Delta$ is a $3\times 3$ diagonal matrix with diagonal:
$$
C(1,0)= \left\{\frac{1}{36} (9 \alpha+6 \beta+\gamma),\frac{1}{36} (9 \alpha+6 \beta+\gamma),\frac{1}{9} (3 \beta+\gamma)\right\}
$$

In the adjoint representation $\w=(1,1)$, $-\Delta$ is an $8\times 8$  diagonal matrix with diagonal:
$$
C(1,1)= \left\{\frac{1}{3} (2 \alpha+\beta),\frac{1}{3} (2 \alpha+\beta),\frac{1}{3} (2 \alpha+\beta),\frac{1}{4} (\alpha+2
   \beta+\gamma),\frac{1}{4} (\alpha+2 \beta+\gamma),\frac{1}{4} (\alpha+2 \beta+\gamma),\frac{1}{4} (\alpha+2 \beta+\gamma), \beta\right\}
$$

More generally, consider the difference $ \Delta -  \beta \Delta_0$, this makes the term $L_4.L_4+L_5.L_5+L_6.L_6+L_7.L_7$ disappear:
\begin{equation*}
\begin{split}
 \Delta -  \beta \Delta_0 = & \frac{-1}{12}\left((\alpha-\beta) \left(L_1.L_1+L_2.L_2+L_3.L_3\right) +(\gamma-\beta) L_8.L_8\right)
\end{split}   
\end{equation*}

From the discussion in sect.~\ref{CasimirGeneral}, we identify $L_1.L_1+L_2.L_2+L_3.L_3$ with the $\SU(2)$ quadratic Casimir, of eigenvalue $4 I (I+1)$ in the irreducible representation of isospin\footnote{In particle physics applications, $I$ is called the isospin, and $Y$  the hypercharge, see our paragraph on notations in sect.~\ref{sec:invariantmetrics}}  $I$ (the highest weight component is  $2I$), and  $L_8.L_8$ with the remaining\footnote{The restriction of $h$ to $\U(2)$ is a bi-invariant metric on this subgroup.} quadratic $\U(1)$ operator, of eigenvalue  $3 Y^2$.
For an  arbitrary representation $\w=(o_1, o_2)$  of $\SU(3)$, and  for a representation that appears in the branching from $\w$ to the subgroup $\U(2)$ (remember, see sect.~\ref{RestrictionSubgroup}, that such a term is characterized by the highest weight $2I$ of $\SU(2)$ and a weight $3 Y$ of $\U(1)$), the eigenvalues  of the modified quadratic Casimir operator  ${\Omega_2^{\U(2)}}=-\Delta$, with $\Delta$ given by (\ref{LaplacianU2}), are therefore:
\begin{equation}
\label{eigenvaluesLaplacian}
 C(o_1,o_2; I, Y) = 
 \beta \,   {C_2(o_1, o_2)} + ( (\alpha - \beta)\,  \frac{1}{3} \,  I (I + 1) + (\gamma -  \beta)\,  \frac{1}{4} \, Y^2)
 \end{equation}
where $C_2$  is the eigenvalue of the Casimir element associated with the Killing form of $\SU(3)$.\\
If $\mu^2\in \RR^+$, then, upon scaling of $h^{-1}$ by $\mu^2$, the rhs of the previous equation gets multiplied by the same factor.

Examples. \\
Consider the fundamental irrep of $\SU(3)$ with highest weight $\w=(1,0)$, of dimension $3$.
The eigenvalue of $-\Delta_0$  given by~(\ref{Casimir2KSU3}) is $4/9$ whereas those of $-\Delta$, given by (\ref{eigenvaluesLaplacian}), for each of the terms appearing in the branching rule obtained at the end of sect.~\ref{RestrictionSubgroup}, are sums of three contributions: 
the irrep $[2]_1$ in the branching of $[3]$ has $I=1/2$ and $Y=1/3$  since $[2]_1 \equiv [2 I +1]_{3Y}$,  and it is such that
$\left\{\beta \,  {{C_2}},\frac{1}{3} I (I+1) (\alpha -\beta ),\frac{1}{4} Y^2 (\gamma -\beta )\right\}= \left\{\frac{4 \beta }{9},\frac{\alpha -\beta }{4},\frac{\gamma -\beta }{36}\right\}$, whose sum is $\frac{1}{36} (9 \alpha +6 \beta +\gamma )$;
in the same way, the three contributions for the irrep $[1]_{-2}$ sum to  $\frac{1}{9} (3 \beta +\gamma )$. Both values were expected since we had to recover the eigenvalues of $-\Delta$ given by the diagonal $3\times 3$ matrix  $C(1,0)$ obtained previously.\\
 For the adjoint representation, whose branching to $\U(2)$ was also given  in  sect.~\ref{RestrictionSubgroup},
 the eigenvalue of $-\Delta_0$ is $1$, and those of $-\Delta$, calculated using (\ref{eigenvaluesLaplacian}), coincide, of course, for each of the four terms of the branching decomposition, with the components of the diagonal $8\times 8$ matrix $C(1,1)$ obtained previously.\\
 Let us finally consider the $\SU(3)$ representation of highest weight $\w=(2,1)$, of dimension $[15]$, whose branching to $\U(2)$ was also considered in sect.~\ref{RestrictionSubgroup}.
The eigenvalue of $-\Delta_0$  given by (\ref{Casimir2KSU3}) is $32/3$ and those of $-\Delta$, given by (\ref{eigenvaluesLaplacian}), read respectively $[4]_1: \frac{5 \alpha }{12}+10 \beta +\frac{\gamma }{4}$,
$[3]_{4}: \frac{2}{3} (\alpha +16 \beta -\gamma )$, $[3]_{-2}: \frac{2 \alpha }{3}+9 \beta +\gamma$, $[2]_1: \frac{1}{36} (-\alpha +376 \beta +9 \gamma )$, $[2]_{-5}: \frac{1}{12} (3 \alpha +50 \beta +75 \gamma)$, $[1]_{-2}: \frac{29 \beta }{3}+\gamma$,
 for the six different terms that appear in the branching rule. The multiplicity, which is $15\times 15$ for a bi-invariant metric (right isometry group $\SU(3)$), becomes $15\times (2I+1)$, where the values of $(2I+1)$ are the consecutive members of the list $\( 4,3,3,2,2,1\)$, when the right isometry group is $\U(2)$. 

\bigskip
Eigenvalues of $\Delta$ for left-invariant metrics with other right isometry groups $K$ can be obtained and discussed along similar lines, this study is left to the reader.

\subsection{The cubic Casimir operator}
From the commutation relations of $Lie(G)$, when $G=\SU(3)$, it is straightforward (but cumbersome) to show that the most general central cubic element of the enveloping algebra of $\mathfrak{su}(3)$ is proportional to $\Omega_3 + x \, \frac{3}{2} \, {\widehat \Omega_2}$, with
{\footnotesize
\begin{equation}
\label{CubicCasimir1}
\begin{split}
\Omega_3 &=
  h_ 1. (\overline {L_ {12}} . L_ {12} + \overline {L_ {45}} .   L_ {45} -2 \overline {L_ {67}} .   L_ {67} )/3 +  h_2. (2\overline {L_ {12}} .   L_ {12} -  \overline {L_ {45}} .   L_ {45} -  \overline {L_ {67}} . L_ {67} )/3
+ \overline {L_ {45}} . L_ {67} .  L_ {12} + \overline {L_ {67}} .\overline {L_ {12}} . L_ {45}+ \\
& \overline {L_ {12}} .   L_ {12} - \overline {L_ {67}} . L_ {67} 
+  \frac {1} {3}\left (\frac {1} {3} ( h_2. h_1. h_1 - h_2. h_2. h_1 )+  \frac {2} {9}( h_1. h_1. h_1  - h_2. h_2. h_ 2) + (h_1. h_1 - h_2. h_2 ) +  (h_ 1 - h_ 2) \right) 
\end{split}   
\end{equation}}
where $L_{ij} = (L_i + i\,  L_j)/2$, $h_1=L_3$, $h_2 = (- L_3 + \sqrt{3} L_8)/2$, and $x$ is an arbitrary real parameter.\\
Setting $x=0$ \ie discarding the quadratic Casimir term ${\widehat \Omega_2}$, we are left with an essentially cubic\footnote{The terms appearing in \ref{CubicCasimir1} (taking $x=0$) are indeed purely cubic when written in terms of the chosen generators $L_{ij}$ and $h_j$, although the same expression, when written in terms of the generators $L_a$, as in \ref{CubicCasimir2}, contains terms linear in $L_3$ and $L_8$.} term that we can write as
{\footnotesize
\begin{equation}
\label{CubicCasimir2}
\begin{split}
\Omega_3 &=
\frac{\left(L_1.L_1+L_2.L_2+L_3.L_3\right).L_8}{4 \sqrt{3}}-\frac{\left(L_4.L_4+L_5.L_5+L_6.L_6+L_7.L_7\right).L_8}{8 \sqrt{3}}+\\
& \frac{1}{8} L_3.\left(L_4.L_4+L_5.L_5-L_6.L_6-L_7.L_7\right)+\frac{1}{4}
   \left(L_1.L_4.L_6+L_1.L_5.L_7-L_2.L_4.L_7+L_2.L_5.L_6\right)-\frac{L_8.L_8.L_8}{12
   \sqrt{3}}+\frac{L_3}{2}-\frac{L_8}{2 \sqrt{3}}
   \end{split}
\end{equation}
}
We can now evaluate this expression  in an irreducible representation of highest weight $\w$, with components $(o_1, o_2)$ in the basis of fundamental weights;  calling $C_3$  the eigenvalue of $\Omega_3$,  one finds   
\begin{equation}
\label{CubicCasimirEigval}
C_3= \frac{1}{27} (o_1 - o_2) (3 + 2 o_1 + o_2) (3 + o_1 + 2 o_2)
\end{equation}
Notice that $C_3$ is equal to $20/27$ in the defining representation, and it vanishes in the adjoint, as in all real representations since it is proportional to $(o_1 - o_2$).

\bigskip
One can play with the idea of relaxing the centrality requirement and look for the most general cubic element commuting with the generators of some Lie subgroup, in particular some right isometry group $K$.
For instance, choosing $K=\U(2)$ we impose the vanishing of Lie derivatives with respect to $L_1, L_2, L_3$ and $L_8$, in which case one finds (again the proof is straightforward) that the most general such cubic element, up to scale, can be written as 
${\Omega_3^{\U(2)}}  + 9\, x \, {\Omega_2^{\U(2)}}$, where $x$ is an arbitrary real parameter, where ${\Omega_2^{\U(2)}} = - \Delta$, with  $\Delta$ the Laplacian given by (\ref{LaplacianU2}), and where the remaining operator,  ${\Omega_3^{\U(2)}}$, is
{\footnotesize
\begin{equation}
\label{cubicsu2u1fundam}
\begin{split}
 {}& \frac{1}{72} \left(6 \sqrt{3} \,A (L_1.L_1.L_8+L_2.L_2.L_8+L_3.L_3.L_8\right)-3
   \sqrt{3} \,B \left(L_4.L_4.L_8+L_5.L_5.L_8+L_6.L_6.L_8+L_7.L_7.L_8\right)-2 \sqrt{3} \,C L_8.L_8.L_8+ \\
   & 9 \,U \left(2 L_1.L_4.L_6+2 L_1.L_5.L_7-2 L_2.L_4.L_7+2  L_2.L_5.L_6+L_3.L_4.L_4+L_3.L_5.L_5-L_3.L_6.L_6-L_3.L_7.L_7+4 L_3\right) -12 \,V \sqrt{3} L_8) 
\end{split}
\end{equation}}
In this expression  $A,B,C,U,V$ denote arbitrary real parameters, and by setting all of them equal to $1$, one recovers the essentially cubic Casimir element $\Omega_3$ given previously.
Using $\SU(3)$ left translations, the element ${\Omega_3^{\U(2)}}$ of the universal enveloping algebra defines a cubic differential operator on the group which is $\SU(3)$ left-invariant by construction,  but also right invariant under $K$.

The interested reader will show that, for an  arbitrary representation $\w=(o_1, o_2)$  of $\SU(3)$, and  for a representation that appears in the branching from $\w$ to the subgroup $\U(2)$ (a representation characterized, as in the previous section, by a pair of integers $(2I, 3Y)$), 
 the eigenvalue of the  operator ${\Omega_3^{\U(2)}}$ is equal to
  \begin{equation}
  \label{cubicsu2u1}
  U * C_3 +  ( U_B\, * {C_2}  + U_{AB}\,  * I(I + 1) + U_{BC}\,  * Y^2 ) *  Y    +  U_V  *  Y
  \end{equation}
 where $C_3$ is given by (\ref{CubicCasimirEigval}), $C_2={\widehat C_2}/6$ is obtained from (\ref{Casimir2KSU3}),  and  $U_B= \frac{3}{2} (U-B)$, $U_{AB}= \frac{1}{2} ( 2 A + B - 3 U)$, $U_{BC}=\frac{1}{8}(3 B - 2 C - U)$,  $U_V =  \frac{1}{2} (U-V)$.
 This is the cubic analog of formula \ref{eigenvaluesLaplacian}.
 
{\scriptsize Example: Choose again the fundamental representation $(1,0)$, so  ${C_2} =4/9$ and $C_3=20/27$.  The branching of this $\SU(3)$ irrep to $\U(2)$ reads $[3]\mapsto [2]_1\oplus[1]_{-2}$ (the notation on the rhs is $[2I+1]_{3Y}$, like in sect.~\ref{RestrictionSubgroup}).
 Equation~\ref{cubicsu2u1fundam}  with $L_a$ replaced by $\lambda_a$ gives a $3\times3$ diagonal matrix with diagonal $\left(
 \frac{A}{4}-\frac{B}{12}-\frac{C}{108}+\frac{3 U}{4}-\frac{V}{6} ,
  \frac{A}{4}-\frac{B}{12}-\frac{C}{108}+\frac{3 U}{4}-\frac{V}{6} ,
 \frac{B}{3}+\frac{2 C}{27}+\frac{V}{3} 
\right)$; 
 its matrix elements can be also be obtained from (\ref{cubicsu2u1}) by setting $I=1/2, Y=1/3$ for the first two, and $I=0, Y=-2/3$ for the last. Finally, one recovers the same value $20/27$ of the undeformed cubic Casimir $C_3$ by setting $A=B=C=U=V=1$.}

\subsection{Sectional curvatures}

Sectional curvatures $\chi$ are associated with the choice of a two-dimensional linear subspace of the tangent space at some point of the manifold under consideration. 
Here the manifold is a Lie group, and the chosen metrics are left-invariant, so it is enough to consider sectional curvatures at the origin.
There is a general formula, due to Milnor \cite{Milnor}, that expresses these quantities in terms of the structure constants of an orthonormal basis (for the chosen metric) of left-invariant vector fields. 
However, we remind the reader that, in these notes, for all choices of the right isometry group $K$ we have chosen the same basis $(X_a)$ to perform curvature calculations, namely the one that is orthonormal for the Killing metric, and which is therefore not orthonormal  for a general left-invariant metric metric $h$. For this reason  we could not use the Milnor formula (with the exception of the Killing metric) and had to rely on the general expression $$\chi(u,v) =  \frac{\langle {\mathcal R}(v,u) u,v\rangle}{\langle u, u\rangle \langle v, v\rangle - \langle u,v\rangle^2}$$ where $u$ and $v$ are two linearly independent vectors at the origin, where the inner product $\langle \,,\, \rangle$ is defined by $h$,  and\footnote{The signs in the definition of  ${\mathcal R}$ are opposite to \cite{Milnor}, this explains the order of arguments in the expression of $\chi(u,v)$.} where ${\mathcal R}(u,v) = [\nabla_u, \nabla_v] - \nabla_{[u,v]}$.
From the calculated Riemann tensors  one can determine the sectional curvatures $\chi(X_a,X_b)=R_{abab}/(h_{aa}h_{bb} - h_{ab}^2)$.  
We list them below, for the parametrizations of $h$ given in (\ref{parametrizations}) that correspond to right isometry groups $K=\SU(3)$,  $\U(2)$ or $\SO(3)$.
For the groups $K = \U(1)\times \U(1)$, $\U(1)_Y$ and $\{e\}$ (the trivial subgroup), we have also explicit results for sectional curvatures in terms of the parameters specifying the metric, but these expressions are too large to be displayed on paper.
For the $K=\U(1)_I$ family, the results are also too large to be displayed but we shall nevertheless give explicit sectional curvatures for a subfamily.
For the special case of the Lorentzian Einstein metric obtained in sect.~\ref{sec: EinsteinMetrics} we only give numerical values (exact values typically involve specific roots of $15^{th}$ degree polynomials).\\
Warning (again): the choice of a pair $X_a,X_b$ with $a\neq b$ determines a two-dimensional linear subspace at the origin but we remind the reader that our vectors $X_a$ are usually not orthonormal for the chosen metric.
By definition $\chi(X_a,X_b)$ is symmetric in $a$ and $b$, and it is not defined (one can set it equal to $0$) for $a=b$. We give tables of $\chi$ for $a < b$ running from  $1$ to $8$. In some cases we also give the Ricci principal curvatures (one can check that their sum is the already given scalar curvature).

\bigskip
\fbox{$K=\SU(3)$}
~
$
\begin{array}{cccccccc}
 . & \frac{\alpha }{12} & \frac{\alpha }{12} & \frac{\alpha }{48} & \frac{\alpha }{48}
   & \frac{\alpha }{48} & \frac{\alpha }{48} & 0 \\
 . & . & \frac{\alpha }{12} & \frac{\alpha }{48} & \frac{\alpha }{48} & \frac{\alpha
   }{48} & \frac{\alpha }{48} & 0 \\
 . & . & . & \frac{\alpha }{48} & \frac{\alpha }{48} & \frac{\alpha }{48} &
   \frac{\alpha }{48} & 0 \\
 . & . & . & . & \frac{\alpha }{12} & \frac{\alpha }{48} & \frac{\alpha }{48} &
   \frac{\alpha }{16} \\
 . & . & . & . & . & \frac{\alpha }{48} & \frac{\alpha }{48} & \frac{\alpha }{16} \\
 . & . & . & . & . & . & \frac{\alpha }{12} & \frac{\alpha }{16} \\
 . & . & . & . & . & . & . & \frac{\alpha }{16} \\
\end{array}
$

All sectional curvatures are non negative, as expected\footnote{Any compact Lie group admits a bi-invariant metric with non-negative sectional curvatures \cite{Milnor} and there is only one bi-invariant metric (up to scale) on  $\SU(3)$.}.
Some of them vanish, also as expected since the $3$-sphere group $\SU(2)$ is the only simply connected Lie group which admits a left invariant metric of strictly positive sectional curvature \cite{Wallach}.
All  Ricci principal curvatures are equal to $\alpha/4$.

\bigskip

\fbox{$K=\U(2)$}
~
$
\begin{array}{cccccccc}
 . & \frac{\alpha }{12} & \frac{\alpha }{12} & \frac{\beta ^2}{48 \alpha } &
   \frac{\beta ^2}{48 \alpha } & \frac{\beta ^2}{48 \alpha } & \frac{\beta ^2}{48
   \alpha } & 0 \\
 . & . & \frac{\alpha }{12} & \frac{\beta ^2}{48 \alpha } & \frac{\beta ^2}{48 \alpha }
   & \frac{\beta ^2}{48 \alpha } & \frac{\beta ^2}{48 \alpha } & 0 \\
 . & . & . & \frac{\beta ^2}{48 \alpha } & \frac{\beta ^2}{48 \alpha } & \frac{\beta
   ^2}{48 \alpha } & \frac{\beta ^2}{48 \alpha } & 0 \\
 . & . & . & . & -\frac{\beta  (9 \alpha  \beta -16 \alpha  \gamma +3 \beta  \gamma
   )}{48 \alpha  \gamma } & \frac{\beta  (4 \alpha -3 \beta )}{48 \alpha } &
   \frac{\beta  (4 \alpha -3 \beta )}{48 \alpha } & \frac{\beta ^2}{16 \gamma } \\
 . & . & . & . & . & \frac{\beta  (4 \alpha -3 \beta )}{48 \alpha } & \frac{\beta  (4
   \alpha -3 \beta )}{48 \alpha } & \frac{\beta ^2}{16 \gamma } \\
 . & . & . & . & . & . & -\frac{\beta  (9 \alpha  \beta -16 \alpha  \gamma +3 \beta 
   \gamma )}{48 \alpha  \gamma } & \frac{\beta ^2}{16 \gamma } \\
 . & . & . & . & . & . & . & \frac{\beta ^2}{16 \gamma } \\
\end{array}
$

Ricci principal curvatures:\\
 $\left\{\frac{\beta ^2}{12 \alpha }+\frac{\alpha }{6},\frac{\beta ^2}{12 \alpha
   }+\frac{\alpha }{6},\frac{\beta ^2}{12 \alpha }+\frac{\alpha }{6},-\frac{\beta ^2}{8
   \alpha }-\frac{\beta ^2}{8 \gamma }+\frac{\beta }{2},-\frac{\beta ^2}{8 \alpha
   }-\frac{\beta ^2}{8 \gamma }+\frac{\beta }{2},-\frac{\beta ^2}{8 \alpha }-\frac{\beta
   ^2}{8 \gamma }+\frac{\beta }{2},-\frac{\beta ^2}{8 \alpha }-\frac{\beta ^2}{8 \gamma
   }+\frac{\beta }{2},\frac{\beta ^2}{4 \gamma }\right\}$.

\bigskip

\fbox{$K=\SO(3)$}
~
$
\begin{array}{cccccccc}
 . & \frac{\alpha ^2}{12 \beta } & \frac{\alpha  (4 \beta -3 \alpha )}{12 \beta } & \frac{\alpha  (4 \beta -3 \alpha )}{48 \beta } &
   \frac{\alpha ^2}{48 \beta } & \frac{\alpha  (4 \beta -3 \alpha )}{48 \beta } & \frac{\alpha ^2}{48 \beta } & 0 \\
 . & . & \frac{\alpha ^2}{12 \beta } & \frac{\alpha ^2}{48 \beta } & \frac{\beta }{48} & \frac{\alpha ^2}{48 \beta } & \frac{\beta }{48} & 0 \\
 . & . & . & \frac{\alpha  (4 \beta -3 \alpha )}{48 \beta } & \frac{\alpha ^2}{48 \beta } & \frac{\alpha  (4 \beta -3 \alpha )}{48 \beta } &
   \frac{\alpha ^2}{48 \beta } & 0 \\
 . & . & . & . & \frac{\alpha ^2}{12 \beta } & \frac{\alpha  (4 \beta -3 \alpha )}{48 \beta } & \frac{\alpha ^2}{48 \beta } & \frac{\alpha  (4
   \beta -3 \alpha )}{16 \beta } \\
 . & . & . & . & . & \frac{\alpha ^2}{48 \beta } & \frac{\beta }{48} & \frac{\alpha ^2}{16 \beta } \\
 . & . & . & . & . & . & \frac{\alpha ^2}{12 \beta } & \frac{\alpha  (4 \beta -3 \alpha )}{16 \beta } \\
 . & . & . & . & . & . & . & \frac{\alpha ^2}{16 \beta } \\
\end{array}
$
\bigskip

In particular, sectional curvatures for the (Jensen) Einstein metric read: 
$
\quad
\begin{array}{cccccccc}
 . & \frac{1 }{132} & \frac{41 1 }{132} & \frac{41 1 }{528} & \frac{1
   }{528} & \frac{41 1 }{528} & \frac{1 }{528} & 0 \\
 . & . & \frac{1 }{132} & \frac{1 }{528} & \frac{11 1 }{48} & \frac{1
   }{528} & \frac{11 1 }{48} & 0 \\
 . & . & . & \frac{41 1 }{528} & \frac{1 }{528} & \frac{41 1 }{528} &
   \frac{1 }{528} & 0 \\
 . & . & . & . & \frac{1 }{132} & \frac{41 1 }{528} & \frac{1 }{528} &
   \frac{41 1 }{176} \\
 . & . & . & . & . & \frac{1 }{528} & \frac{11 1 }{48} & \frac{1 }{176} \\
 . & . & . & . & . & . & \frac{1 }{132} & \frac{41 1 }{176} \\
 . & . & . & . & . & . & . & \frac{1 }{176} \\
\end{array}
$

Ricci principal curvatures:
$\left\{-\frac{\alpha  (\alpha -2 \beta )}{4 \beta },-\frac{\alpha  (\alpha -2 \beta )}{4
   \beta },-\frac{\alpha  (\alpha -2 \beta )}{4 \beta },-\frac{\alpha  (\alpha -2 \beta
   )}{4 \beta },-\frac{\alpha  (\alpha -2 \beta )}{4 \beta },\frac{5 \alpha ^2+\beta
   ^2}{24 \beta },\frac{5 \alpha ^2+\beta ^2}{24 \beta },\frac{5 \alpha ^2+\beta ^2}{24
   \beta }\right\}$.
   
   For the Jensen metric they are all equal to $21\alpha/44$.

\bigskip

\fbox{$K=\U(1)_I$.} 

The general case, with its $8$ parameters, is too large to be displayed. 
We can nevertheless (using tiny fonts (!)) exhibit the sectional curvatures obtained for the subfamily $\theta=0, \zeta = 0$, $\delta = \gamma$.

{\tiny
\begin{equation*}
\begin{split}
&
 \left\{
 \begin{split}
 & \frac{
    \cdot , \alpha  (4 \beta -3 \alpha )}{12 \beta },\frac{\alpha ^2}{12 \beta
   },\frac{\gamma ^5-2 \eta ^2 \gamma ^3-4 \alpha  \eta ^2 \gamma ^2+\eta ^4 \gamma +4
   \alpha ^2 \eta ^2 \gamma +4 \alpha  \eta ^4}{48 \alpha  \gamma ^3-48 \alpha  \gamma 
   \eta ^2},\frac{\gamma ^5-2 \eta ^2 \gamma ^3-4 \alpha  \eta ^2 \gamma ^2+\eta ^4
   \gamma +4 \alpha ^2 \eta ^2 \gamma +4 \alpha  \eta ^4}{48 \alpha  \gamma ^3-48 \alpha 
   \gamma  \eta ^2},\frac{\gamma ^5-2 \eta ^2 \gamma ^3-4 \alpha  \eta ^2 \gamma ^2+\eta
   ^4 \gamma +4 \alpha ^2 \eta ^2 \gamma +4 \alpha  \eta ^4}{48 \alpha  \gamma ^3-48
   \alpha  \gamma  \eta ^2},  \\
& \qquad  \frac{\gamma ^5-2 \eta ^2 \gamma ^3-4 \alpha  \eta ^2 \gamma
   ^2+\eta ^4 \gamma +4 \alpha ^2 \eta ^2 \gamma +4 \alpha  \eta ^4}{48 \alpha  \gamma
   ^3-48 \alpha  \gamma  \eta ^2},0 
  \end{split}
  \right\},
  \\
&   \left\{
 \begin{split}
&   \cdot ,    \cdot  , \frac{\alpha ^2}{12 \beta },
    \frac{\gamma ^5-2 \eta ^2 \gamma ^3-4 \alpha  \eta ^2 \gamma ^2+\eta ^4 \gamma +4
   \alpha ^2 \eta ^2 \gamma +4 \alpha  \eta ^4}{48 \alpha  \gamma ^3-48 \alpha  \gamma 
   \eta ^2},\frac{\gamma ^5-2 \eta ^2 \gamma ^3-4 \alpha  \eta ^2 \gamma ^2+\eta ^4
   \gamma +4 \alpha ^2 \eta ^2 \gamma +4 \alpha  \eta ^4}{48 \alpha  \gamma ^3-48 \alpha 
   \gamma  \eta ^2},\frac{\gamma ^5-2 \eta ^2 \gamma ^3-4 \alpha  \eta ^2 \gamma ^2+\eta
   ^4 \gamma +4 \alpha ^2 \eta ^2 \gamma +4 \alpha  \eta ^4}{48 \alpha  \gamma ^3-48
   \alpha  \gamma  \eta ^2},\\
   & \frac{\gamma ^5-2 \eta ^2 \gamma ^3-4 \alpha  \eta ^2 \gamma
   ^2+\eta ^4 \gamma +4 \alpha ^2 \eta ^2 \gamma +4 \alpha  \eta ^4}{48 \alpha  \gamma
   ^3-48 \alpha  \gamma  \eta ^2},0
    \end{split}
   \right\},\\
   & \left\{
      \cdot ,    \cdot ,    \cdot , \frac{\gamma ^2-\eta ^2}{48 \beta
   },\frac{\gamma ^2-\eta ^2}{48 \beta },\frac{\gamma ^2-\eta ^2}{48 \beta },\frac{\gamma
   ^2-\eta ^2}{48 \beta },0\right\},\\
 &  \left\{
  \begin{split}
  &   \cdot ,    \cdot ,    \cdot ,    \cdot , 
 \frac{\beta  \left(-9 \gamma ^4+16 \epsilon 
   \gamma ^3+18 \eta ^2 \gamma ^2-16 \epsilon  \eta ^2 \gamma -9 \eta ^4+8 \alpha 
   \epsilon  \eta ^2\right)-3 \epsilon  \left(\gamma ^2-\eta ^2\right)^2}{48 \beta 
   \gamma ^2 \epsilon },\frac{-4 \alpha ^2 \eta ^2-3 \left(\gamma ^2-\eta ^2\right)^2+4
   \alpha  \left(\gamma ^3-\eta ^2 \gamma +3 \epsilon  \eta ^2\right)}{48 \alpha  \gamma
   ^2},\frac{4 \alpha ^2 \eta ^2-3 \left(\gamma ^2-\eta ^2\right)^2+4 \alpha 
   \left(\gamma ^3-\gamma  \eta ^2\right)}{48 \alpha  \left(\gamma ^2-\eta
   ^2\right)},\\
   &
    \frac{\gamma ^5-2 \eta ^2 \gamma ^3-4 \epsilon  \eta ^2 \gamma ^2+\eta ^4
   \gamma +4 \epsilon ^2 \eta ^2 \gamma +4 \epsilon  \eta ^4}{16 \gamma ^3 \epsilon -16
   \gamma  \epsilon  \eta ^2}
     \end{split}\right\},\\
   &
   \left\{
      \cdot ,    \cdot ,    \cdot ,    \cdot ,    \cdot , \frac{4 \alpha ^2 \eta ^2-3 \left(\gamma
   ^2-\eta ^2\right)^2+4 \alpha  \left(\gamma ^3-\gamma  \eta ^2\right)}{48 \alpha 
   \left(\gamma ^2-\eta ^2\right)},
   \frac{-4 \alpha ^2 \eta ^2-3 \left(\gamma ^2-\eta
   ^2\right)^2+4 \alpha  \left(\gamma ^3-\eta ^2 \gamma +3 \epsilon  \eta ^2\right)}{48
   \alpha  \gamma ^2},\frac{\gamma ^5-2 \eta ^2 \gamma ^3-4 \epsilon  \eta ^2 \gamma
   ^2+\eta ^4 \gamma +4 \epsilon ^2 \eta ^2 \gamma +4 \epsilon  \eta ^4}{16 \gamma ^3
   \epsilon -16 \gamma  \epsilon  \eta ^2}\right\},\\
   &
   \left\{
      \cdot ,    \cdot ,    \cdot ,    \cdot ,    \cdot ,    \cdot , \frac{\beta  \left(-9 \gamma
   ^4+16 \epsilon  \gamma ^3+18 \eta ^2 \gamma ^2-16 \epsilon  \eta ^2 \gamma -9 \eta
   ^4+8 \alpha  \epsilon  \eta ^2\right)-3 \epsilon  \left(\gamma ^2-\eta ^2\right)^2}{48
   \beta  \gamma ^2 \epsilon },\frac{\gamma ^5-2 \left(\gamma ^2+2 \epsilon  \gamma -2
   \epsilon ^2\right) \eta ^2 \gamma +(\gamma +4 \epsilon ) \eta ^4}{16 \gamma  \epsilon 
   \left(\gamma ^2-\eta ^2\right)}\right\},\\
   & 
   \left\{
   \cdot ,    \cdot ,    \cdot ,    \cdot ,    \cdot ,    \cdot ,    \cdot , 
   \frac{\gamma ^5-2 \left(\gamma ^2+2
   \epsilon  \gamma -2 \epsilon ^2\right) \eta ^2 \gamma +(\gamma +4 \epsilon ) \eta
   ^4}{16 \gamma  \epsilon  \left(\gamma ^2-\eta ^2\right)}\right\}
   \end{split}
\end{equation*}
}

The scalar curvature, for this family, was given by (\ref{scalarcurvatureU1I}).

\bigskip

\fbox{$K=\U(1)_I$.  Special case: Lorentzian Einstein metric}

{\scriptsize
$
\begin{array}{cccccccc}
 . & 0.156539 & 0.0589315 & 0.0207884 & 0.0207884 & 0.0207884 & 0.0207884 & 0 \\
 . & . & 0.0589315 & 0.0207884 & 0.0207884 & 0.0207884 & 0.0207884 & 0 \\
 . & . & . & 0.000619797 & 0.000619797 & 0.000619797 & 0.000619797 & 0 \\
 . & . & . & . & 0.108778 & -0.0335267 & 0.0410926 & -0.054309 \\
 . & . & . & . & . & 0.0410926 & -0.0335267 & -0.054309 \\
 . & . & . & . & . & . & 0.108778 & -0.054309 \\
 . & . & . & . & . & . & . & -0.054309 \\
\end{array}
$
}

Other features of this metric have been discussed in sect.\ref{sec: EinsteinMetrics}.

\subsection{Ricci decomposition  (examples)}
The Ricci decomposition of the Riemann tensor associated with the Levi-Civita connection defined by the metric $h$, 
namely  $R = C + \frac{1}{d-2} (\rho - \frac{\tau}{d} h) \circleland h + \frac{\tau}{2d(d-1)} h \circleland h$, where $d$ is the dimension (here $d=8$),  $C$ is the Weyl tensor, and $\circleland$ denotes the Kulkarni-Nomizu product of two $(0,2)$ tensors,
expresses the Riemann tensor  (here thought of as a $(0,4)$ tensor),  as an orthogonal direct sum. 
Such a decomposition can be considered for an arbitrary metric, in particular for homogeneous metrics.
As a verification of our calculations involving curvatures, we have checked this identity for all the familes of metrics considered in this paper.  It would be of course paper-consuming to list all the non-zero entries of the relevant tensors, 
nevertheless, in a few cases, it may be useful to mention a numerical consequence of this identity, namely the following norm decomposition: 
\begin{equation}
|R|^2 = |C|^2 + |\frac{1}{d-2} (\rho - \frac{\tau}{d} h) \circleland h|^2 + |\frac{\tau}{2d(d-1)} h \circleland h|^2
\end{equation}
In those cases where the results are reasonably short we give $|R|^2$ followed by a triple containing the three  contributions, in the same order as  in the previous equation.

$K=\SU(3)$, (bi-invariant metrics):
$\left\{\frac{\alpha ^2}{2},\{\frac{5 \alpha ^2}{14},0,\frac{\alpha^2}{7}\}\right\}$
\bigskip

$K=\U(2)$:
{\scriptsize
\begin{equation*}
\begin{split}
&\{\frac{8 \alpha ^4 \gamma ^2+\alpha ^2 \beta ^2 \left(51 \beta ^2-144 \beta  \gamma +176 \gamma ^2\right)+6 \alpha  \beta ^3 \gamma  (5 \beta -16 \gamma )+23 \beta ^4 \gamma ^2}{96 \alpha ^2 \gamma ^2},\\
&\{\frac{80 \alpha ^4 \gamma ^2-24 \alpha ^3 \beta  \gamma  (\beta -8 \gamma )+\alpha ^2 \beta ^2 \left(909 \beta ^2-2448 \beta  \gamma +2600 \gamma ^2\right)+6 \alpha  \beta ^3 \gamma  (79 \beta -240 \gamma )+377 \beta ^4 \gamma ^2}{2016 \alpha ^2 \gamma ^2},\\
&\frac{20 \alpha ^4 \gamma ^2+12 \alpha ^3 \beta  \gamma  (\beta -8 \gamma )+\alpha ^2 \beta ^2 \left(45 \beta ^2-144 \beta  \gamma +236 \gamma ^2\right)+6 \alpha  \beta ^3 \gamma  (7 \beta -24 \gamma )+29 \beta ^4 \gamma ^2}{576 \alpha ^2 \gamma ^2},
\frac{1}{448} \left(-\frac{\beta ^2}{\alpha }+2 \alpha +\beta  \left(8-\frac{\beta }{\gamma }\right)\right)^2\}\}
\end{split}
\end{equation*} 
}     
      
$K=\SO(3)$:
$\left\{\frac{295 \alpha ^4-660 \alpha ^3 \beta +460 \alpha ^2 \beta ^2+\beta ^4}{192 \beta ^2},
\{
\frac{5 \left(508 \alpha ^4-1110 \alpha ^3 \beta +733 \alpha ^2 \beta ^2+12 \alpha  \beta ^3+\beta ^4\right)}{2016 \beta ^2},
\frac{5 \left(11 \alpha ^2-12 \alpha  \beta +\beta ^2\right)^2}{2304 \beta ^2},
\frac{\left(-5 \alpha ^2+20 \alpha  \beta +\beta ^2\right)^2}{1792 \beta ^2}\}\right\}$

In particular, for the Einstein solution (Jensen case: $\beta = 11 \alpha$), we have  $\left\{\frac{2639 \alpha ^2}{968},\{\frac{2135 \alpha ^2}{968},0,\frac{63 \alpha ^2}{121}\}\right\}$.
\bigskip

$K=\U(1)_I$: The results giving the norms are too large to be displayed. We only consider the Lorentzian Einstein solution. In that case, one can obtain these four norms as roots of appropriate 15th degree polynomials with (very) large integer coefficients. We shall not print them. 
Numerically, one finds $\{0.115257, \{0.0813543, 0, 0.0339023\}\}$.
\bigskip

The expressions of the four norms, for $K=\U(1)\times \U(1)$ and for $K=\U(1)_Y$, are also very large, and we shall not display them.
\bigskip

Notice that for all Einstein spaces the square norm $|(\rho - \frac{\tau}{d} h) \circleland h|^2 $ vanishes, as it should.


\section{Physical applications}
\subsection{Particle physics  and the GMO formula}
\paragraph{Physical considerations.}
In the standard model of elementary particles, more precisely in the quark model based on the Lie group $G=\SU(N)$,  the particles called mesons are associated with the space of intertwiners ($G$-equivariant morphisms) from $V_f \otimes {\overline V_f}$, where $V_f$ is the defining representation and ${\overline V_f}$ its conjugate, to the irreducible representations that appear in the product, namely the adjoint representation, and the trivial one. This decomposition, for $\SU(3)$, in terms of dimensions, read $[3]\otimes[\overline 3] = [8] \oplus [1]$.
Mesons are described as basis vectors of the tensor product, specified by appropriate $G$-Clebsch-Gordan coefficients (or by the  $3J$ Wigner symbols of the group $G$) associated with the chosen spaces of intertwiners, but this is not our concern here.

Quarks (resp. anti-quarks) are basis vectors in the space of the defining representation of $G$ (resp. its conjugate), and mesons are called ``bound states '' of a quark and an anti-quark.
In particle physics parlance $G$ is the ``flavor group'', which, in the presently accepted model, means $\SU(N)$, and where $N$ can only be $2,3,4,5$ or $6$. 
When $N=2$, $G$ is called the isospin group and the basis vectors of the defining representation (the quarks) are nicknamed `up' and `down'. 
When $N=3$ (in this paper we restrict our attention to $G=\SU(3)$) they are nicknamed `up', `down', `strange',  and are respectively denoted by $u,d,s$.

The classical (\ie not quantum field theoretical) description of mesons, as sections of appropriate vector bundles over the space-time manifold,  also specifies their behavior with respect to space-time symmetries (\ie under action of the Lorentz group or of the Poincar\'e group), but we don't have to be more precise here, it is enough to say that there are several families of mesons differing by their space-time properties, the two most important families being the so-called pseudo-scalar mesons and the vector mesons. 

Quarks are not observable but mesons are, and they have masses. Experimentally the pseudo-scalar mesons have masses that are close (same remark for the vector mesons), and this is precisely the reason why, historically, they were described as members of the same Lie group multiplet (basis vectors of some irrep).
Calculating meson masses in terms of more fundamental parameters is a task that goes beyond the possibilities of (perturbative) quantum field theory, in particular of quantum chromodynamics, 
but it remains that, phenomenologically, one can assume that interactions of mesons,  in particular the quadratic operator responsible for their masses, or their mass splitting, commutes with the generators of the chosen flavor group, or of a subgroup of the latter. This hypothesis is at the origin of several mass relations.

Experimentally, particles of the same $\SU(3)$ multiplet have approximately the same masses, but this is even more so when they are members of the same irreducible representation of the $\U(2)$ subgroup (locally $\SU(2) \times \U(1)$) defined in table (\ref{specificK}).
For this reason, it is natural to describe (or approximate) the unknown mass operator as the Laplacian on $\SU(3)$ associated with an appropriate left-invariant metric for which the right isometry group is $K=\U(2)$, and to look at the consequences of this ansatz.
We rescale the dual metric  by $\mu^2$ to fix the dimensions ($\mu$ will have the dimensions of a mass).
The eigenvalues of the Laplacian are given by (\ref{eigenvaluesLaplacian}), we write them again below (the whole expression is now multiplied by $\mu^2$). To an irrep  $(o_1,o_2)$ of $\SU(3)$ branching to irreps of $\U(2)$ labelled by the isospin value $I$ and hypercharge $Y$ we associate the square mass 
$m^2= \mu ^2 \, C(o_1,o_2)$ given by

\begin{equation}
m^2  =   {{C_2}(o_1, o_2)} \, \beta \, \mu ^2 + \left(\frac{1}{3} I (I+1)  (\alpha-\beta) \mu ^2+\frac{1}{4} Y^2 (\gamma-\beta) \mu ^2\right)
\end{equation}

\paragraph{Pseudo-scalar mesons and the Gell-Mann-Okubo formula.}

The branching of the adjoint representation (octet) of $\SU(3)$, when restricted to the previously defined $\U(2)$ subgroup, in terms of $\SU(2)$ highest weights, reads\footnote{the components are written in the basis of fundamental weights.}: $(1,1) \rightarrow (2)+(1)+(1)+(0)$; this is often written in terms of dimensions  of irreps ($[8] \rightarrow [3]_0 + [2]_3 + [2]_{-3} + [1]_0$), with the same notations $[2I+1]_{3Y}$ as in sect.~\ref{RestrictionSubgroup}, since there is no possible confusion in the present case. 
The corresponding mesons are the three pions $\{\pi^+, \pi^0, \pi^-\}$, for which $I=1$, $Y=0$, the four kaons
$\{K^+, K^0\}$ for which $I=1$, $Y=1$,  and $\{\overline{K^0}, K^-\}$ for which $I=1$, $Y=-1$, and the eta particle, for which $I=0$, $Y=0$.

In the adjoint representation of $\SU(3)$,  ${C_2}=1$, so that one obtains immediately:\
$$
\left\{m^2{}_{\pi }, \, m^2{}_K,\, m^2{}_{\eta }\right\}=\left\{\frac{1}{3}  (2 \alpha+\beta)\mu ^2, \frac{1}{4} (\alpha+2 \beta+\gamma)\mu ^2, \beta \mu ^2\right\}
$$
Experimentally:  $m_ {\pi^+} = m_ {\pi^-} = 139.57 \, \text {MeV}$,
    $m_ {\pi^0} = 134.976\, \text {MeV}$,
     $m_ {K^+} = m_ {K^-} = 493.677 \, \text {MeV}$,
      $m_ {K^0} =  m_ {\overline {K^0}} = 497.64 \, \text {MeV}$, and
      $m_\eta = 549\,  \text{MeV}$. For pions and kaons we use averaged masses $m_\pi \simeq 137\,  \text{MeV}$,  $m_K \simeq 496\,  \text{MeV}$, $m_\eta = 549\,  \text{MeV}$.
The corresponding values\footnote{These values where already obtained in \cite{RCGEF}, up to a scaling factor equal to $12$ coming from the fact that the bi-invariant metric used in that reference for normalizing purposes was a multiple of the Killing metric.} of parameters are then 
$\mu^2 \alpha_{exp} \simeq - (350\, \text{MeV})^2,\, \mu^2 \beta_{exp} \simeq (549\, \text{MeV})^2, \, \mu^2 \gamma_{exp} \simeq (710\, \text{MeV})^2$. 
Their ratios, normalized by $\beta$,  are $(\alpha_{exp},\beta_{exp},\gamma_{exp})/\beta_{exp} = (-0.406591, 1, 1.67156)$.

We stress the fact that the above is nothing else than an educated fit: it is neither a prediction nor a ``post-diction'' since the number of unknown parameters is the same as the number of values coming from experiment.
In order to get a prediction, one needs at least one more relation between the parameters $\alpha, \beta, \gamma$; such a relation (expressed in a rather different way) was postulated in the sixties, by making the hypothesis that the $\SU(3)$ the mass operator could be well approximated by keeping only its singlet and octet components, therefore neglecting the contribution from the representation of dimension $27$ (see for instance \cite{Swart}). In our language, this amounts to neglect the  ${}_{27}h^{-1}$ component of the (dual) pseudo\footnote{The signature of the bilinear form, with the previous values of 
$\alpha_{exp}, \beta_{exp} , \gamma_{exp}$,  is $(5,3)$.}-Riemannian metric in the decomposition $h^{-1} = {{}_1h^{-1}}  + {{}_8h^{-1}} + {{}_{27}h^{-1}}$ 
discussed in sect.~\ref{metricdecomposition}.
In other words the coefficient $C=\frac{1}{40} (\alpha -4 \beta +3 \gamma)$ given in (\ref{decomposition1827}) is set to $0$. One can then eliminate the parameter $\gamma$, for example, and obtain
$$\left\{m^2{}_{\pi }, \, m^2{}_K,\, m^2{}_{\eta }\right\} \simeq \left(\frac{1}{3}  (2 \alpha+\beta) \mu ^2,\, \frac{1}{6}  (\alpha+5 \beta)\mu ^2,\, \beta \mu ^2 \right),$$
$$\text{which implies the relation:} {\hskip 1.cm} m_{\eta }^2 \simeq \frac{1}{3} \left(4 m_K^2-m_{\pi }^2\right).$$ This is the celebrated Gell-Mann-Okubo formula for pseudo-scalar mesons (the formula using square masses, see for instance the article on GMO formulae in Wikipedia). It holds reasonably well\footnote{The first published mass relation of this type (1961) was for baryons, and in particular for hyperons of the decuplet. As it is well known this equation lead to the discovery of the $\Omega^-$ particle (and to the Nobel Prize in Physics 1969). In our language, this latter formula, which is linear in masses (not quadratic)  could be obtained from the eigenvalues of a Dirac operator on $\SU(3)$ associated with a left-invariant metric for which the right isometry group is $\U(2)$.}, although, using the experimental values of $m_{\pi}$ and $m_{K}$, it  leads to a value of $567\, \text{MeV}$ for the mass of the $\eta$, which is slightly too big.

\smallskip

{\sl Remarks}: Rather than using a (rough) Laplacian, for some appropriate $\SU(3)\times \U(2)$ invariant metric, one could be tempted of using the Yamabe (conformal) Laplacian that differs\footnote{The conformal Laplacian is $\Delta - \tfrac{d-2}{4(d-1)}\, \tau$, where $\tau$ is the scalar curvature. In our case $d=8$.}  from the latter by a simple shift (equal to  $ - \tfrac{3}{14}\,  \frac{1}{4} \left(-\frac{\beta ^2}{\alpha }+2 \alpha +\beta  \left(8-\frac{\beta }{\gamma }\right)\right)$).
but this does not seem to lead to anything physically particularly interesting. 
One could  also play with the idea of using analogous considerations to study other particle multiplets, to generalize the previous analysis to Lie groups $SU(N)$ for $N>3$, or even to consider other kinds of ``symmetry breaking" scenarios (selecting other right-invariant isometry groups), but this would lie beyond the intended scope of these notes.

\paragraph{Warning.} The physicist reader certainly knows, and the mathematician reader should be warned, that the above way of obtaining mass relations for some elementary particles, from considerations on Laplacians (or Dirac operators) associated with left-invariant metrics (actually $\SU(3) \times \U(2)$ invariant metrics on the Lie group $\SU(3)$ is not standard, in the sense that, although not a new observation (see \cite{RCGEF}),  this approach is not widely known and it is not the way it is taught. It can be noticed that, whatever the starting point one chooses (the elementary quark model, or more sophisticated approaches like QCD or chiral perturbation theory),  the elementary mathematical considerations leading to these mass relations are similar: they involve representation theory of $\SU(3)$, the branching to $\SU(2)\times \U(1)$ (of isospin and hypercharge), the fact that masses should be related to eigenvalues of linear or quadratic  operators (would-be Hamiltonian or mass operators) that are not explicitly known, and an approximation of ``octet dominance'' (in the present case it is the hypothesis that one can neglect, in the metric,
 a contribution associated with the $27$ dimensional representation of $\SU(3)$  ---see sect~\ref{metricdecomposition} and the discussion in \cite{Swart}).  
 However, one thing is to cook up a physico-mathematical formula that works at least approximately in some particular cases, another is to derive it in the framework of a physical theory. The contemporary attitude is to view GMO mass formulae for hadrons as remote consequences of a fundamental theory called ``The Standard Model''; this theory, in its usual formulation, does not explicitly involve  considerations on the left-invariant geometry of Lie groups $\SU(N)$ ---this $N$ standing for the number of ``flavors''. Our observation that the same mass formulae can be interpreted as expressions describing the spectrum of appropriate differential operators associated with particular left-invariant metrics on Lie groups (it is not difficult to extend the previous results to $N>3$) may not be significant but, 
notwithstanding this possibility, 
the result suggests that it could be, or should be, justified, as a formal consequence of some currently accepted physical theory (a subject that we don't investigate in these notes) and maybe trigger some new developments of the latter.

\subsection{Other possibly physical considerations}
The fact that some particular left-invariant metrics, which may be Riemannian or pseudo-Riemannian, can sometimes be Einstein metrics played no role in the previous discussion on particle masses.
Now, for the last thirty years or so, many theoretical physicists have been concerned with the construction of classical field models, quantum field theories, and string theories, incorporating, on top of space and time,  several  ``extra-dimensions'' that we do not perceive (the prototype being the old Kaluza-Klein theory).
Such models, that are often speculative, are described by equations that sometimes require the total space of the theory (or maybe a quotient of the latter) to be an Einstein manifold.
We have no wish to comment this endeavor and shall refrain from suggesting anything in that direction but one may notice that, if needed, 
the examples of eight dimensional Einstein manifolds described in sect~\ref{sec: EinsteinMetrics} can be used to construct higher dimensional pseudo-Riemannian spaces, of given signature, that are also Einstein.
For instance one can build a Lorentzian homogeneous Einstein metric on the $11$-dimensional compact manifold $\SU(3)\times \SU(2)$
where the first factor of the Cartesian product is endowed with the Lorentzian Einstein structure described previously and where the second factor (diffeomorphic with the sphere $S^3$) is endowed with its standard $\SO(4)$ invariant Einstein metric.

 \end{document}